\documentclass[twoside,reqno]{amsart} 

\usepackage{array}
\usepackage{times}            %
\usepackage{amsthm}           %
\usepackage{amssymb}          %
\usepackage{enumerate}        %
\usepackage{euscript}         %
\usepackage{dan2e}
\usepackage{varioref}         %

\usepackage{ifpdf}

\usepackage[matrix,arrow,curve,frame]{xy}      %
\SilentMatrices
\ifpdf
\else
  \xyoption{ps}
\fi
\xymatrixcolsep{1.7pc} 
\xymatrixrowsep{1.7pc}
\newdir{ >}{{}*!/-5pt/\dir{>}}               %
\newdir{ |}{{}*!/-5pt/\dir{|}}                  %
\newcommand{\circar}{\ar@{}|<>(.42)*\cir<1.7pt>{}} %

\usepackage{hyperref}

\numberwithin{equation}{section}

\newenvironment{roenumerate}{\begin{enumerate}[\ulp i\urp]}{\end{enumerate}}

\newtheorem{thm}[equation]{Theorem} %
\newtheorem{defn}[equation]{Definition}
\newtheorem{prop}[equation]{Proposition}
\newtheorem{cor}[equation]{Corollary}
\newtheorem{lemma}[equation]{Lemma}
\newtheorem{cnj}[equation]{Conjecture}

\theoremstyle{definition}  %
\newtheorem{example}[equation]{Example}
\newtheorem{question}[equation]{Question}
\newtheorem{note}[equation]{Note}
\newcommand{\ulp}{\textup{(}}
\newcommand{\urp}{\textup{)}}

\newcommand{\dfn}{\textbf} %
\newcommand{\mdfn}[1]{\dfn{\mathversion{bold}#1}} %

\newcommand{\cat}               {\EuScript}  
\renewcommand{\cA}                {{\cat A}}

\renewcommand{\cC}                {{\cat C}}
\renewcommand{\cD}                {{\cat D}}
\renewcommand{\cE}                {{\cat E}}
\renewcommand{\cF}                {{\cat F}}

\renewcommand{\cI}                {{\cat I}}
\newcommand{\cJ}                {{\cat J}}
\newcommand{\cK}                  {{\cat K}}

\renewcommand{\cP}                {{\cat P}}
\newcommand{\cQ}                {{\cat Q}}
\newcommand{\cR}                {{\cat R}}
\renewcommand{\cS}                {{\cat S}}
\newcommand{\cT}                {{\cat T}}
\newcommand{\cU}                {{\cat U}}

\newcommand{\Ch}                {\textup{Ch}}

\newcommand{\K}                 {\cK}

\newcommand{\D}                 {\cD}

\DeclareMathOperator{\pdim}              {proj.dim}

\DeclareMathOperator{\Stl}       {St}

\DeclareMathOperator{\RP}        {\R P}
\newcommand{\tighten}{\mskip-0.4\thinmuskip}
\newcommand{\RPn}       {{\RP^{\tighten n}}}

\newcommand{\RPi}       {{\RP^{\tighten\infty}}}
\newcommand{\RPnp}      {{\RP^{\tighten n}_{\tighten 0}}}

\newcommand{\RPip}      {{\RP^{\tighten\infty}_{\tighten 0}}}

\DeclareMathOperator{\CP}        {\C P}

\newcommand{\CPi}       {{\CP^{\tighten\infty}}}

\newcommand{\CPip}      {{\CP^{\tighten\infty}_{\tighten 0}}}

\newcommand{\Quat}{\field H}
\DeclareMathOperator{\HP}        {\Quat P}
\newcommand{\Oct}{\field O}
\DeclareMathOperator{\OP}        {\Oct P}

\newcommand{\EM}                {Eilenberg--Mac\,Lane\xspace}

\newcommand{\Zp}{\Z \mspace{-1.0mu}/\mspace{-1.0mu} p}
\newcommand{\Zt}{\Z \mspace{-1.0mu}/\mspace{-1.0mu} 2}
\newcommand{\Hp}{H\Zp}
\newcommand{\Ht}{H\Zt}

\newcommand{\nul}               {\textup{-null}}
\newcommand{\proj}              {\textup{-proj}}
\newcommand{\exact}             {\textup{-exact}}
\newcommand{\epi}               {\textup{-epi}}

\newcommand{\om}                {\infty}

\newcommand{\Proj}              {\cP}
\newcommand{\Projnull}          {\Proj\nul}

\newcommand{\Projexact}         {\Proj\exact}
\newcommand{\Projexactproj}     {\Proj\exact\proj}
\newcommand{\Projepi}           {\Proj\epi}

\newcommand{\Projo}             {{\Proj_{\!\om}}}
\newcommand{\I}                 {\cI} 
\newcommand{\Iproj}             {\I\proj}

\newcommand{\Io}                {{\I^\om}}
\newcommand{\Exact}             {\cC}
\newcommand{\Exactproj}         {\Exact\proj}
\newcommand{\Exactprojexact}    {\Exact\proj\exact}

\newcommand{\Epi}               {\cE} 
\newcommand{\Epiproj}           {\Epi\proj}

\newcommand{\PC}[1]{\cP \cC(#1)}
\newcommand{\PCS}{\PC{\cS}}

\newcommand{\zero}{\mathbf{0}}
\newcommand{\one}{\mathbf{1}}

\newcommand{\LEQ}{\leq}
\newcommand{\J}{\cJ}
\newcommand{\Qroj}{\cQ}
\newcommand{\Rroj}{\cR}

\DeclareMathOperator{\Pure}     {P}
\renewcommand{\PExt}              {\Pure\!\Ext}
\newcommand{\PzExt}             {\Pure^0\!\Ext}
\newcommand{\PoExt}             {\Pure^1\!\Ext}

\newcommand{\PkExt}             {\Pure^k\!\Ext}
\newcommand{\PlExt}             {\Pure^l\!\Ext}
\newcommand{\PE}[1]{\Pure^{#1}\!\Ext}
\newcommand{\RMod}              {\textup{$R$-Mod}\xspace}
\newcommand{\kay}               {k}

\newcommand{\Ai}{$A_{\infty}$\xspace}

\newcommand{\invlim}            {\varprojlim\nolimits}
\newcommand{\invlimone}         {{\invlim^{\!1}}}

\newcommand{\invlimi}           {\invlim^{\!i}}
\newcommand{\invlims}           {\invlim^{\!s}}
\newcommand{\dirlim}            {\varinjlim}

\newcommand{\floor}[1]          {\lfloor #1 \rfloor}

\DeclareMathOperator{\obj}{obj}
\DeclareMathOperator{\mor}{mor}

\newcommand{\PBS}[1]{\let\blick=\\#1\let\\=\blick}

\newcommand{\period}    {{\makebox[0pt][l]{\hspace{2pt} .}}}
\newcommand{\comma}     {{\makebox[0pt][l]{\hspace{2pt} ,}}}

\newcommand{\mysection}{\section}

\begin{document}

\title{Ideals in triangulated categories: phantoms, ghosts and skeleta}
\author{J.\ Daniel Christensen}
\thanks{The author was partially supported by an NSF grant and
  an NSERC scholarship.}
\thanks{This version is more up to date than the published version.  
See the last page for details.}
\address{Department of Mathematics\\ University of Western Ontario\\
London, Ontario, Canada}
\email{jdc@uwo.ca}



\vspace*{-.4in}

\begin{abstract}
  We begin by showing that 
  in a triangulated category, specifying a projective class is equivalent 
  to specifying an ideal $\I$ of morphisms with certain properties, and
  that if $\I$ has these properties, then so does each of its powers.  
  We show how a projective class leads to an Adams spectral sequence
  and give some results on the convergence and collapsing of this
  spectral sequence.
  We use this to study various ideals.  In the stable homotopy category
  we examine phantom maps, skeletal phantom maps, superphantom maps,
  and ghosts.  
  (A ghost is a map which induces the zero map of homotopy groups.)
  We show that ghosts lead to a stable analogue of the Lusternik--Schnirelmann
  category of a space, and we calculate this stable analogue for 
  low-dimensional real projective spaces.
  We also give a relation between ghosts and
  the Hopf and Kervaire invariant problems.
  In the case of \Ai modules over an \Ai ring spectrum, the ghost spectral
  sequence is a universal coefficient spectral sequence.
  From the phantom projective class we derive
  a generalized Milnor sequence for filtered
  diagrams of finite spectra, and from this it follows that the group of
  phantom maps from $X$ to $Y$ can always be described as a $\invlimone$
  group.
  The last two sections focus on algebraic examples.
  In the derived category of an abelian category 
  we study the ideal of maps inducing the zero map of homology groups
  and find a natural setting for a result of Kelly on the vanishing
  of composites of such maps.
  We also explain how pure exact sequences relate to phantom maps in the
  derived category of a ring and give an example showing that phantoms
  can compose non-trivially.
\end{abstract}

\date{February 25, 2013}

\maketitle

\tableofcontents

\pagestyle{plain}

\renewcommand{\baselinestretch}{1.40}\normalsize

\mysection{Introduction}\label{se:intro}

Let $\cS$ be a triangulated category, such as the stable homotopy
category or the derived category of a ring.
One often tries to study an object $X$ in $\cS$ by building it
up from a class of better understood pieces.
When this is done, there may be maps $X \ra Y$ that aren't seen by the
pieces (this is made precise below).
Obvious questions arise, such as how efficiently $X$ can be built
from the given class of pieces, and how the unseen maps behave
under composition.
This paper presents a systematic way of studying such phenomena,
namely by showing that they are captured in the notion of a ``projective
class''.  We then apply this formalism to various interesting examples.

If $\Proj$ is a collection of objects of $\cS$, denote by $\Projnull$
the collection of all maps $X \ra Y$ such that the composite
$P \ra X \ra Y$ is zero for all objects $P$ in $\Proj$ and all
maps $P \ra X$.
These are the maps that the objects of $\Proj$ fail to ``see''.
If $\I$ is a collection of maps of $\cS$, denote by $\Iproj$ 
the collection of all objects $P$ such that the composite
$P \ra X \ra Y$ is zero for all maps $X \ra Y$ in $\I$ and all
maps $P \ra X$.
A \dfn{projective class} is a pair $(\Proj,\I)$ with
$\Projnull = \I$ and $\Iproj = \Proj$ such that for each
object $X$ there is a cofibre sequence $P \ra X \ra Y$ with 
$P$ in $\Proj$ and with $X \ra Y$ in $\I$.
The objects of $\Proj$ are referred to as \dfn{projectives}.

This definition of ``projective class'' appears on the surface
to be different from other definitions that have appeared,
but in the next section we show that it is in fact equivalent.
However, when working in a triangulated category, we claim
that the above definition is more natural than the others.
For example, the collection $\I$ of maps is automatically
an \dfn{ideal} in $\cS$.  That is,
if $f$ and $g$ are parallel maps in $\I$, then $f + g$ is in $\I$.
And if $f$, $g$ and $h$ are composable and $g$ is in $\I$,
then both $f g$ and $g h$ are in $\I$.  
(All of our ideals will be two-sided.)
Many commonly occurring ideals in fact form part of a projective class.

Our definition also has the feature that it leads to a sequence of
``derived'' projective classes.
The powers $\I^n$ of the ideal $\I$ form a decreasing
filtration of the class of maps of $\cS$.
Let $\Proj_1 = \Proj$ and inductively define
$\Proj_n$ to be the class of all retracts of objects
$Y$ that sit in a cofibre sequence $X \ra Y \ra P$
with $X$ in $\Proj_{n-1}$ and $P$ in $\Proj$.
In this way we get an increasing filtration of the class of objects
of $\cS$, and we have the following result.

\begin{thm}\label{th:pc}
  For each $n$, the pair $(\Proj_n,\I^n)$ is a projective class.
\end{thm}

This result is a special case of the ``product'' operation
on projective classes.  
Both are discussed in more detail in Section~\ref{se:operations}.

Theorem~\ref{th:pc} forms the cornerstone of our investigations of various
ideals in the stable homotopy category and the derived category
of an abelian category, 
investigations which occupy us in the later sections of the paper.  
But before getting to these examples, we show in Section~\ref{se:ass} how 
a projective class leads to an Adams spectral sequence, and we prove that
this spectral sequence is conditionally convergent~\cite{bo:ccss}
if the ideal is closed under countable coproducts
and the projective class is \dfn{generating}, \ie if for each non-zero
$X$ there is a projective $P$ and a non-zero map $P \ra X$.
The \dfn{length} of an object $X$ is the smallest $n$ such that 
$X$ is in $\Proj_n$, and is infinite if $X$ is not in any $\Proj_n$.
If $X$ has finite length and we have a generating projective class
such that $\I$ is closed under countable coproducts,
then the Adams spectral sequence abutting to $[X,Y]$ is strongly
convergent~\cite{bo:ccss}.
Our next observation is that
when dealing with a generating projective class,
an upper bound on the length of an object $X$ is given by 
$1 + \pdim X$, where the projective dimension of $X$ is 
the length of the shortest projective resolution of $X$
with respect to the given projective class.

The remaining sections deal with examples of projective classes.
Section~\ref{se:phantoms} studies phantom maps in any category
which has a generating phantom projective class (Definition~\ref{de:gppc}).
A phantom map is a map which is zero when restricted to any
``finite'' object (Definition~\ref{de:finite}).
A necessary tool is the introduction of weaker variants of colimits,
such as weak colimits and minimal cones.
We conclude that if our category satisfies ``Brown Representability''
(Definition~\ref{de:br}), then every object has projective dimension
at most one, the composite of any two phantoms is zero, and the
Adams spectral sequence degenerates into a generalized Milnor sequence.

In Section~\ref{se:top-phantoms} we examine various types of phantom
maps in the stable homotopy category.
We begin with ordinary phantom maps, drawing on the abstract
work done in the previous section.
As a special case of the generalized Milnor sequence, we
get the following result.
\begin{thm}\label{th:1.2}
  Let $X$ be a CW-spectrum and let $\{X_\alpha\}$ be the filtered diagram
  of finite CW-subspectra of $X$.
  For any spectrum $Y$ there is a short exact sequence
  \[ 0 \lra \invlimone [\Sigma X_\alpha, Y] \lra [X,Y] 
       \lra \invlim [X_\alpha, Y] \lra 0 . \]
  The kernel consists precisely of the phantom maps.
  Moreover, $\invlimi [\Sigma X_\alpha, Y]$ vanishes for $i \geq 2$.
\end{thm}
This result is straightforward when $X$ has finite skeleta,
but is more delicate in general.  It was also proved by
Tetsusuke Ohkawa~\cite{oh:vtaybkss}.
For further historical comments, see the paragraph following
Theorem~\ref{th:milnor}.

In the next part of Section~\ref{se:top-phantoms} 
we study skeletal phantom maps, maps which
are zero when restricted to any ``skeleton'' of the source
(Definition~\ref{de:cellular-tower}).
And in the final part
we prove that ``superphantoms'' exist, answering a question of
Margolis~\cite[p.~81]{ma:ssa}.  
A superphantom is a map which is zero when restricted to any 
(possibly desuspended) suspension spectrum.

We turn to the study of maps which induce the zero map in homotopy
groups in Section~\ref{se:ghosts}.
(We dub these maps ``ghosts''.) 
Here things are much more interesting in that our ideal has infinite order.
If a spectrum $X$ has length $k$ with respect to this ideal,
then any composite of $k$ Steenrod operations in the mod 2 cohomology
of $X$ is zero.
For this reason, we view the length of a spectrum with respect to the
ideal of ghosts
as a stable analogue of the Lusternik-Schnirelmann
category of a space.
In the second part of the section we focus on calculating 
the lengths of real projective spaces.
We give upper and lower bounds on the length of $\RPn$ which agree
for $2 \leq n \leq 19$, and 
we show that the filtration of the Steenrod squares is closely
related to the Hopf and Kervaire invariant problems.
In the third and final part of the section, we show that
the Adams spectral sequence with respect to the ghost projective
class in the category of \Ai modules over an \Ai ring spectrum
is a universal coefficient spectral sequence, and we show how
this gives another theoretical lower bound on the length of a
spectrum.

The last two sections deal with the algebraic analogues of
phantoms and ghosts.
In the derived category $\D$ of an abelian category, a \dfn{ghost} is
a map which induces zero in homology, and the ideal of ghosts is
a natural ideal to study.
We show that if our abelian category has enough projectives and
satisfies Grothendieck's AB 5 axiom,
then the ideal of ghosts is part of a generating projective class, 
and we prove the following result.
\begin{thm}\label{th:1.3}
  Let $X$ be a complex such that the projective dimensions of
  $B_n X$ and $H_n X$ are less than $k$ for each $n$.
  Then the projective dimension of $X$ 
  with respect to the ideal of ghosts
  is less than $k$.
  In particular, $X$ has length at most $k$, and a $k$-fold composite
  \[ X \lra Y^1 \lra \cdots \lra Y^k \]
  of maps each zero in homology is zero in $\D$.
\end{thm}
One can deduce from this a result of Max Kelly~\cite{ke:cmizhm},
whose work provided the inspiration for the general framework
presented in this paper.  See Section~\ref{se:alg-ghosts} for details.

In the last section we discuss phantom maps in the derived category
of an associative ring.  A map $X \ra Y$
is phantom if and only if the composite $W \ra X \ra Y$ is
zero (in the derived category) for each bounded complex $W$ of finitely
generated projectives.
Our main result here is a relation between phantom
maps and pure extensions of $R$-modules as well as an analog
of Theorem~\ref{th:1.3}.
We reproduce an example which shows that when
$R = \C[x,y]$, Brown representability doesn't hold in the
derived category.

I have had the pleasure of discussing this work with many mathematicians,
to whom I owe a great debt.  
I mention in particular 
Mike Hopkins, for help with the result that $\I^n$ forms part of
a projective class;
Neil Strickland, for the joint work described in Section~\ref{se:phantoms},
which led to the view of a projective class described in this paper;
Mark Mahowald, for the efficient construction of $\RPi$;
and my advisor, Haynes Miller, 
for his constant encouragement, support and good advice.

\mysection{Projective classes} \label{se:projclass}

This section has three parts.  In the first part we recall
the definition of a projective class given by Eilenberg and
Moore~\cite{eimo:frha}.  This definition focuses on the relation between
``projective'' objects and three-term ``exact'' sequences.
It appears to be the most general notion
in that it allows one to do rudimentary homological algebra
in any pointed category.
In the second part we show that the more familiar
relation between projectives and ``epimorphisms'' can
be used to define ``projective class'' as long as the
pointed category in which we are working has weak kernels.
Finally, in the third part, we show that it is equivalent
to use the relation between projectives and ``null maps''
as long as our pointed category has weak kernels and in
addition has the property that every map \emph{is} a weak kernel.
Any triangulated category satisfies these conditions.

A reader who is more interested in the examples may move on
to the next section, using the definition of ``projective class'' 
from the introduction.

Everything we say can be dualized to give a discussion
of injective classes.

\subsection[In pointed categories]{Projective classes in pointed categories} \label{sse:pointed}

We recall the notion of a projective class in 
a pointed category.
We will be brief; the elegant original paper~\cite{eimo:frha} leaves no room
for improvement.

A category $\cS$ is \dfn{pointed} if it contains an object
which is both initial and terminal.  If such an object exists,
it is unique up to isomorphism and is denoted $0$.
For any two objects $X$ and $Y$ of $\cS$, there is a
unique map from $X$ to $Y$ that factors through $0$, and
we denote this map by $0$ as well.
By using the zero map as a basepoint, the hom functor
$\cS(-,-)$ takes values in the category of pointed sets.

For the rest of this section we assume that $\cS$ is pointed.

A composable pair of maps $X \ra Y \ra Z$ is said to be a 
\dfn{(length two) complex} if the composite is zero.

A complex $X \llra{f} Y \llra{g} Z$ in the category of pointed sets is
\dfn{exact} if $f(X) = g^{-1}(*)$, where $*$ denotes the basepoint in $Z$.

\begin{defn}
  Let $\Proj$ be a collection of objects of $\cS$.
  A complex $X \ra Y \ra Z$ such that
  \[ \cS(P,X) \lra \cS(P,Y) \lra \cS(P,Z) \]
  is an exact sequence of pointed sets for each $P$ in $\Proj$
  is said to be \mdfn{$\Proj$-exact}, and
  the collection of all such complexes is denoted $\Projexact$.

  Now let $\Exact$ be any collection of length two complexes.
  An object $P$ such that
  \[ \cS(P,X) \lra \cS(P,Y) \lra \cS(P,Z) \]
  is an exact sequence of pointed sets for each complex 
  $X \ra Y \ra Z$ in $\Exact$ is said to be \mdfn{$\Exact$-projective},
  and the collection of all such objects is denoted $\Exactproj$.
\end{defn}

Note that the class of $\Exact$-projectives is closed under coproducts
and retracts.

\begin{defn}
  Let $\Proj$ be a collection of objects in 
  $\cS$ and let $\Exact$ be a collection of length two complexes.
  The pair $(\Proj,\Exact)$ is \dfn{complementary} if
  $\Projexact = \Exact$ and $\Exactproj = \Proj$,
  and is a \dfn{projective class} if, in addition,
  for each morphism $X \ra Y$ in $\cS$
  there is a morphism $P \ra X$ such that $P$ is in $\Proj$ and
  $P \ra X \ra Y$ is in $\Exact$.
\end{defn}

It is easily checked that $(\Projexactproj,\Projexact)$ 
is a complementary pair for any collection $\Proj$ of objects of $\cS$.
Similarly, $(\Exactproj,\Exactprojexact)$ is complementary for any collection
$\Exact$ of length two complexes.

If $\cS$ is an additive category, then all of the usual results
about projective resolutions may be proved.
See~\cite{eimo:frha} or~\cite{humo:dghasv} for details.

\subsection[In categories with weak kernels]{Projective classes in categories with weak kernels} \label{sse:wk}

We show that in a category with weak kernels 
we can define ``projective class'' by using epimorphisms instead
of exact sequences.
This observation was made in~\cite{eimo:frha} for a category with
strict kernels.

Given a map $f : X \ra Y$ in a pointed category $\cS$, 
a \dfn{weak kernel} for $f$ is a map $W \ra X$ such that 
\[ \cS(V,W) \lra \cS(V,X) \lra \cS(V,Y) \]
is an exact sequence of pointed sets for each $V$ in $\cS$.
This says that a map $V \ra X$ factors through $W$ if and only if the 
composite $V \ra X \ra Y$ is zero.
In particular, the composite $W \ra X \ra Y$ is zero.
We say that $\cS$ \dfn{has weak kernels}
if every map in $\cS$ has a weak kernel.

\begin{defn}
  Let $\Proj$ be a collection of objects of $\cS$.
  A map $X \ra Y$ such that
  \[ \cS(P,X) \lra \cS(P,Y) \]
  is an epimorphism for each $P$ in $\Proj$
  is said to be \mdfn{$\Proj$-epic}, and
  the collection of all such maps is denoted $\Projepi$.

  Now let $\Epi$ be any collection of maps of $\cS$.
  An object $P$ such that
  \[ \cS(P,X) \lra \cS(P,Y) \]
  is an epimorphism for each map $X \ra Y$ in $\Epi$
  is said to be \mdfn{$\Epi$-projective},
  and the collection of all such objects is denoted $\Epiproj$.
\end{defn}

Note that the class of $\Epi$-projectives is closed under coproducts
and retracts.

\begin{prop} \label{pr:epi}
  Let $\cS$ be a pointed category with weak kernels.
  Let $\Proj$ be a collection of objects and $\Epi$ a collection
  of morphisms such that $\Projepi = \Epi$ and $\Epiproj = \Proj$.
  Assume also that for each $X$ there is a map $P \ra X$ in $\Epi$ with
  $P \in \Proj$.
  Then $(\Proj,\Projexact)$ is a projective class.
  Moreover, every projective class is of this form for a unique pair 
  $(\Proj,\Epi)$ satisfying the above conditions.
\end{prop}

\begin{proof}
  Suppose that we are given a pair $(\Proj,\Epi)$ with
  $\Projepi = \Epi$ and $\Epiproj = \Proj$ such that for each
  $X$ there is a map $P \ra X$ in $\Epi$ with $P$ in $\Proj$.
  Let $\Exact = \Projexact$.
  First we will prove that $(\Proj,\Exact)$ is complementary.
  By definition $\Exact = \Projexact$, and it
  is clear that $\Proj \subseteq \Exactproj$, so all that remains
  to be shown is that $\Exactproj \subseteq \Proj$.
  Let $X$ be a $\Exact$-projective object, and choose a map
  $P \ra X$ in $\Epi$ with $P$ in $\Proj$.
  The sequence $P \ra X \ra 0$ is $\Proj$-exact and so it is
  exact under $\cS(X,-)$.  Thus $X$ is a retract of $P$ and therefore
  is in $\Proj$.
 
  To finish the proof of the first part of the proposition, we must show that
  there are enough projectives.
  Let $f : Y \ra Z$ be a map and choose a weak kernel $X \ra Y$ for $f$.
  Let $P \ra X$ be a map in $\Epi$ with $P$ in $\Proj$.
  Then it is clear that $P \ra Y \ra Z$ is $\Proj$-exact, and therefore
  $(\Proj,\Projexact)$ is a projective class.

  Now we prove the converse.  Suppose $(\Proj,\Exact)$ is a projective
  class.  If this projective class is obtained from a pair $(\Proj,\Epi)$
  as above, then we have $\Epi = \Projepi$, so uniqueness is clear.
  Thus our task is to show that taking this as a definition of $\Epi$
  we have that $(\Proj,\Epi)$ satisfies the hypotheses of the first
  part of the proposition.
  By definition $\Epi = \Projepi$, and it is clear that
  $\Proj \subseteq \Epiproj$, so we must show that 
  $\Epiproj \subseteq \Proj$.
  Let $X$ be in $\Epiproj$.  Choose a map $P \ra X$ 
  with $P$ in $\Proj$ so that $P \ra X \ra 0$ is in $\Exact$.
  It is easy to see that $P \ra X$ is in $\Epi$, and since $X$ is
  in $\Epiproj$, $X$ is a retract of $P$.
  Therefore $X$ is in $\Proj$.

  All that is left to be done is to show that for each $X$ there is a 
  map $P \ra X$ in $\Epi$ with $P$ in $\Proj$.
  As above, simply choose a map $P \ra X$ so that $P \ra X \ra 0$
  is in $\Exact$.
\end{proof}

\subsection[In triangulated categories]{Projective classes in triangulated categories}

In this section we assume that our pointed category $\cS$ has 
weak kernels and in addition has the property that every map
is a weak kernel.  Another way to say this is that every
map $X \ra Y$ lies in a sequence
\[ W \ra X \ra Y \ra Z \]
which is exact under $\cS(U,-)$ for all objects $U$.
A triangulated category satisfies this condition with
$W$ the fibre of $X \ra Y$ and $Z$ the cofibre
(so $Z \iso \Sigma W$).
There are various references for triangulated categories.
The reader already looking at Margolis' book~\cite{ma:ssa} will find 
Appendix~2 to be a handy reference.
A standard (and good) reference is Verdier's portion of
SGA $4 \frac{1}{2}$~\cite{ve:cd}.

\begin{defn}
  Let $\Proj$ be any collection of objects of $\cS$.
  A map $X \ra Y$ such that
  \[ \cS(P,X) \lra \cS(P,Y) \]
  is the zero map for each $P$ in $\Proj$
  is said to be \mdfn{$\Proj$-null}, and
  the collection of all such maps is denoted $\Projnull$.

  Now let $\I$ be any collection of maps of $\cS$.
  An object $P$ such that
  \[ \cS(P,X) \lra \cS(P,Y) \]
  is the zero map for each map $X \ra Y$ in $\I$
  is said to be \mdfn{$\I$-projective},
  and the collection of all such objects is denoted $\Iproj$.
\end{defn}

Note that the class of $\I$-projectives is closed under coproducts 
and retracts, and that 
if $f$, $g$ and $h$ are composable and $g$ is in $\Proj$-null, 
then both $f g$ and $g h$ are in $\Proj$-null.
If the category $\cS$ is additive, then $\Proj$-null is an
ideal.

\begin{prop}
  Let $\cS$ be a pointed category with weak kernels such that
  every map is a weak kernel.
  Let $\Proj$ be a collection of objects and $\I$ a collection
  of morphisms such that $\Projnull = \I$ and $\Iproj = \Proj$.
  Assume also that for each $X$ there is a projective $P$ and a 
  map $P \ra X$ which is a weak kernel of a map in $\I$.
  Then $(\Proj,\Projexact)$ is a projective class.
  Moreover, every projective class is of this form for a unique pair 
  $(\Proj,\I)$ satisfying the above conditions.
\end{prop}

\begin{proof}
  This can be proved directly, paralleling the proof of the
  analogous result in the previous part of this section, but it is simpler
  and more illuminating to show how this formulation relates to the
  epimorphism formulation and then to apply Proposition~\ref{pr:epi}.

  Let $\Proj$ be a class of objects, and let $X \ra Y \ra Z$
  be a sequence exact under $\cS(U,-)$ for each $U$ in $\Proj$.
  Then $X \ra Y$ is $\Proj$-epic if and only if $Y \ra Z$ is
  $\Proj$-null.
  Using this, and the facts that every $\Proj$-epic map can
  be detected in this way (since every map \emph{is} a weak kernel) and
  that every $\Proj$-null map can be detected in this way
  (since every map \emph{has} a weak kernel), it is
  easy to see that pairs $(\Proj,\I)$ as described in
  the hypotheses correspond bijectively to pairs $(\Proj,\Epi)$
  as described in the previous part of this section.
\end{proof}

For the rest of the paper we will be working in a triangulated
category and will freely make use of the equivalent ways of
thinking of a projective class.
Also, in a triangulated category we can replace the condition
\begin{quote}
  to any $X$ we can associate a projective $P$ and a map $P \ra X$ which 
  is a weak kernel of a map in $\I$
\end{quote}
with the condition
\begin{quote}
  any $X$ lies in a cofibre sequence
  $P \ra X \ra Y$
  with $P$ in $\Proj$ and $X \ra Y$ in $\I$.
\end{quote}
Indeed, the latter clearly implies the former.
And given a weak kernel $P \ra X$ of a map in $\I$, it is easy to
check that the cofibre $X \ra Y$ of $P \ra X$ is in $\I$.

\mysection{Operations on projective classes}
\label{se:operations}

For the rest of this paper, $\cS$ will be a triangulated category
containing all set-indexed coproducts.
We will sometimes slip and call a coproduct a ``wedge'', and
we will write \mbox{$X \Wedge Y$} for the coproduct of $X$ and $Y$.
All of our projective classes will be \dfn{stable}. That is, both
$\Proj$ and $\I$ are assumed to be closed under suspension and
desuspension.  If $(\Proj,\I)$ is a projective class in $\cS$, we will
call the objects of $\Proj$ \dfn{projective}.

\subsection{Meets and products}

There is a natural ordering on the class $\PCS$ of projective classes in $\cS$.
For projective classes $(\Proj,\I)$ and $(\Qroj,\J)$,
write $(\Proj , \I) \LEQ (\Qroj ,\J)$ if $\I $ is contained in $\J $.
The projective class $\zero = (\obj \cS , 0 )$, 
whose ideal contains only the zero maps, is the smallest projective class.
The projective class $\one = (0,\mor \cS)$, 
whose ideal contains all maps, is the largest projective class.

\begin{prop}\label{pr:meets}
  Let $\{(\Proj_{\alpha },\I_{\alpha})\}$ be a set of projective classes.
  Then the intersection
\[
\bigcap_{\alpha } \I_{\alpha}
\]
  is an ideal which forms part of a projective class.
  The projectives are precisely the retracts of wedges of objects
  from the union
\[
\bigcup_{\alpha } \Proj_{\alpha} .
\]
\end{prop}

Note that not every ideal forms part of a projective class,
so there is some content to this proposition.

The following lemma will be used to prove the proposition.  In
fact, it will be used just about every time we prove that we
have an example of a projective class.

\begin{lemma} \label{le:pc}
  Let $\Proj$ be a class of objects closed under retracts and
  let $\I$ be an ideal.  Assume that $\Proj$ and $\I$ are
  orthogonal, \ie that the composite $P \ra X \ra Y$ is zero
  for each $P$ in $\Proj$, each map $X \ra Y$ in $\I$, and
  each map $P \ra X$.
  Also assume that each object $X$ lies in a cofibre sequence
  $P \ra X \ra Y$ with $P$ in $\Proj$ and $X \ra Y$ in $\I$.
  Then $(\Proj,\I)$ is a projective class.
\end{lemma}

\begin{proof}[Proof of Lemma]
  All that we have to show is that $\I\proj \subseteq \Proj$
  and that $\Proj\nul \subseteq \I$.
  For the former, assume that $X$ is in $\I\proj$ and choose
  a cofibre sequence $P \ra X \ra Y$ with $P$ in $\Proj$ and
  $X \ra Y$ in $\I$.  Since $X$ is in $\I\proj$, the map
  $X \ra Y$ is zero, and so $X$ is a retract of $P$.  Hence
  $X$ is in $\Proj$, as we have assumed that the latter is
  closed under retracts.

  To show that $\Proj\nul \subseteq \I$ is equally easy, using
  that $\I$ is an ideal.
\end{proof}

\begin{proof}[Proof of Proposition~\ref{pr:meets}]
  Let $\I$ denote the intersection of the ideals $\I_{\alpha}$ and
  let $\Proj$ denote the collection of retracts of wedges of objects each
  of which lies in some $\Proj_{\alpha}$.
  It is clear that $\Proj$ and $\I $ are orthogonal, so we
  must verify that each object $X$ lies in a cofibre sequence
  $P \ra X \ra Y$ with $P$ in $\Proj$ and $X \ra Y$ in $\I$.
  Let $X$ be an object of $\cS$.  
  For each $\alpha$ choose a cofibre sequence 
  $P_{\alpha} \ra X \ra Y_{\alpha }$.
  Consider the map $\Wedge P_{\alpha } \ra X$.
  The cofibre of this map is zero when restricted to each $P_{\alpha }$
  and so must lie in each $\I_{\alpha }$.
  Thus we have a cofibre sequence
  $P \ra X \ra Y$ with $P$ in $\Proj$ and $X \ra Y$ in $\I$.
\end{proof}

The projective class $(\Proj ,\I )$ constructed in the
proposition is the \dfn{meet} of the
set $\{(\Proj_{\alpha },\I_{\alpha })\}$.
That is, it is the greatest lower bound.

I don't know whether an arbitrary set of projective classes has a \dfn{join}
(\ie a least upper bound).  If $A$ is a set of projective classes and
the collection of all projective classes which are upper bounds for
$A$ has a meet, then this is the join of $A$.  However, the collection
of all projective classes which are upper bounds for $A$ might not be
a set, and so might not have a meet.  (The proof above involved a
coproduct over all $\alpha$.)

There is also a product on $\PCS$.

\begin{prop}\label{pr:products}
If the ideals $\I $ and $\J $ are parts of projective classes $(\Proj,\I)$ and
$(\Qroj,\J)$, then so is their
product $\I \J $, which consists of all composites $f g$
with $f$ in $\I $ and $g$ in $\J $.
The projectives are precisely those objects which are retracts
of objects $X$ which lie in cofibre sequences
$Q \ra X \ra P$ with $Q$ in $\Qroj $ and $P$ in $\Proj $.
\end{prop}

\begin{proof}
Write $\Rroj$ for the collection of retracts
of objects $X$ which lie in cofibre sequences
$Q \ra X \ra P$ with $Q$ in $\Qroj $ and $P$ in $\Proj $.
Then $\Rroj$ is closed under retracts and coproducts.
Using the fact that in an additive category
binary coproducts are biproducts, one can show that $\I \J $ is in
fact an ideal.  (The point is that it is automatically closed
under sums of parallel maps.)
We will show now that $\Rroj$ and $\I \J $ are orthogonal, \ie that
any composite
\[
W \lra X \llra{f} Y \llra{g} Z
\]
is zero if
$W$ is in $\Rroj $, $f$ is in $\J $, and $g$ is in $\I $.
We can assume without loss of generality that $W$ lies in 
a cofibre sequence $Q \ra W \ra P$ with $Q$ in $\Qroj $ and
$P$ in $\Proj $.
In the following diagram
\[ \xymatrix{
   Q \ar[d] \\
   W \ar[d] \ar[r] & X \ar[r]^{f} & Y \ar[r]^{g} & Z \\
   P \ar@{{}-->}[rru]
} \]
  the dashed arrow exists because $\Qroj$ and $\J$ are
  orthogonal.  
  The map $P \ra Y \ra Z$ is zero because
  $\Proj$ and $\I$ are orthogonal, so the map $W \ra X \ra Y \ra Z$
  is zero.
  This shows that $\Rroj$ and $\I \J$ are orthogonal.

  It remains to show that any $X$ lies in a cofibre sequence
  $W \ra X \ra Z$ with $W$ in $\Rroj $ and $X \ra Z$ in $\I \J$.
  To do this, choose a cofibre sequence
  $Q \ra X \ra Y$ with $Q$ in $\Qroj $ and $X \ra Y$ in $\J $.
  Now choose a cofibre sequence
  $P \ra Y \ra Z$ with $P$ in $\Proj $ and $Y \ra Z$ in $\I $,
  giving a diagram
\[
  \xymatrixcolsep{0.6pc} 
  \xymatrixrowsep{1.05pc}
  \begin{array}{c}  %
  \xymatrix{ \relax
   X \ar[rr] && Y \ar[ld] \circar[ld] \ar[rr] && Z \ar[ld] \circar[ld] \\
   & Q \ar[lu] && P \period \ar[lu]
  }
  \end{array}
\]
  (A circle on an arrow $A \ra B$ denotes a map $A \ra \Sigma B$.)
  Let $W$ be the fibre of the composite $X \ra Y \ra Z$.
  Using the octahedral axiom, one sees that $W$ lies in a cofibre
  sequence $Q \ra W \ra P$.
  Thus $W \ra X \ra Z$ is the sequence we seek.

  With the help of the lemma, we have proved that
  $(\Rroj,\I \J )$ is a projective class.
\end{proof}

Here are some formal properties of the intersection and product operations, 
all of which are easy to prove.

\begin{prop}\label{pr:prod}
Let $(\Proj,\I)$, $(\Qroj,\J)$ and $(\Rroj,\K)$ be projective classes.  Then
\begin{roenumerate}
\item $\zero \I = \zero = \I \zero$ and 
$\zero \cap \I = \zero = \I \cap \zero$.  
\ulp Recall that $\zero = (\obj \cS , 0 )$.\urp
\item $\one \I = \I = \I \one$ and $\one \cap \I = \I = \I \cap \one $.  
\ulp Recall that $\one = (0,\mor \cS)$.\urp
\item $\I \cap \I = \I $ and $\I \cap \J = \J \cap \I $.
\item $\I \J \LEQ \I \cap \J$ and $\J \I \LEQ \I \cap \J$.
\item If $\J \LEQ \K $, then $\I \J \LEQ \I \K$, $\J \I \LEQ \K \I$,
$\I \cap \J \LEQ \I \cap \K $ and $\J \cap \I \LEQ \K \cap \I $.      \qed
\end{roenumerate}
\end{prop}

To make the proposition readable, we have blurred the distinction
between a projective class and the ideal that it corresponds to.
We hope that by now the reader has been convinced that
the ideal is the most natural part of a projective class.
Indeed, it is usually easier to describe the ideal than the projectives,
and it is through the ideal that the operations on projective
classes arise naturally.

We next describe the two filtrations that a projective class determines.

\subsection{Two filtrations}

Fix a stable projective class $(\Proj , \I )$.
Define $\I^n$ to be the collection of all $n$-fold composites
of maps in $\I$.  
The ideals $\I^n$ form a decreasing filtration
of the class of morphisms of $\cS$;
write $\Io$ for the intersection.
We will use the notation $\I(X,Y)$ for \mbox{$\I \cap \cS(X,Y)$}, and
more generally $\I^n(X,Y)$ for $\I^n \cap \cS(X,Y)$, $1 \leq n \leq \om$.

By the results of the previous part of this section, each of these
ideals forms part of a projective class.  To fix notation and
terminology, we will explicitly describe the increasing filtration 
of the class of objects.
Let $\Proj_1 = \Proj$ and inductively define
$\Proj_n$ to be the class of all retracts of objects
$Y$ which sit in cofibre sequences $X \ra Y \ra P$
with $X$ in $\Proj_{n-1}$ and $P$ projective.
If $X$ is in $\Proj_n$ but not in $\Proj_{n-1}$ we say that
$X$ has \dfn{length} $n$ with respect to the projective class $(\Proj,\I)$.
We also say that $X$ can be \dfn{built} from $n$ objects of $\Proj$.
Write $\Projo$ for the class of all retracts of wedges of
objects of finite length.  
We write $\Proj_0$ for the collection of
zero objects of $\cS$, and say that they have length $0$.
For symmetry, we write $I^0$ for the collection of all morphisms
in $\cS$.

The next theorem follows immediately from Propositions~\ref{pr:meets} 
and~\ref{pr:products}.

\begin{thm} \label{th:powers}
  For $0 \leq n \leq \om$, the pair $(\Proj_n,\I^n)$ is a projective 
  class.                                                            \qed
\end{thm}

We say that the projective classes $(\Proj_{n},\I^{n})$ are
\dfn{derived} from $(\Proj,\I)$. 

\begin{note}
  If there is a cofibre sequence
  \[ X \ra Y \ra Z \]
  with $X$ of length $m$ and $Z$ of length $n$,
  then $Y$ has length at most $m+n$.
  For if we have a composite
  \[ Y \lra Y_1 \lra \cdots \lra Y_{m+n} \]
  of $m+n$ maps each in $\I$, the composite of the first $m$ maps
  is zero when restricted to $X$, and so factors through $Z$.
  But then the composite $Y \ra Z \ra Y_m \ra \cdots \ra Y_{m+n}$
  is zero because $Z$ has length $n$.
  So, by the theorem, $Y$ is in $\Proj_{m+n}$.
\end{note}

\begin{note}
  Suppose $\I$ and $\J $ are projective classes with $\I \LEQ \J $.
  Since products respect order (Proposition~\ref{pr:prod}), 
  we have that $\I^n \LEQ \J^n$.
  Therefore the length of an object $X$ with respect to $\I $ is
  no more than the length of $X$ with respect to $\J $.
  In contrast, the \dfn{projective dimension} of an object $X$
  (the length of the shortest projective resolution of $X$)
  might not respect the order.
\end{note}

Much of this paper will focus on studying these filtrations,
both abstractly and in particular examples.

\mysection{The Adams spectral sequence}\label{se:ass}

Associated to a projective class $(\Proj,\I)$ in a triangulated category $\cS$
is an Adams spectral sequence which we now describe.
The Adams spectral sequence was discussed in the same generality
in~\cite{mi:aaanss}.

Let $X$ be an object.  By repeatedly using the fact that
there are enough projectives, one can form a diagram
\begin{equation} \label{eq:ar} 
  \quad
  \xymatrixcolsep{0.6pc} 
  \xymatrixrowsep{1.05pc}
  \begin{array}{c}  %
  \xymatrix{ \relax
   \makebox[2em][r]{$X = X_0$}
       \ar[rr] && X_1 \ar[ld] \circar[ld] \ar[rr] && X_2 \ar[ld] \circar[ld] \ar[rr] && X_{3} \ar[ld] \circar[ld] \\
   & P_0 \ar[lu] && P_1 \ar[lu] && P_2 \ar[lu] 
  }
  \end{array}
  \cdots
\end{equation}
with each $P_n$ projective, each map $X_n \ra X_{n+1}$ in $\I$,
and each triangle exact.
We call such a diagram an \dfn{Adams resolution} of $X$ with
respect to the projective class $(\Proj,\I)$.
Let $W_n$ be the fibre of the composite map $X \ra X_n$.
Then $W_0 = 0$, $W_1 = P_0$ and $W_n$ sits in a cofibre
sequence $W_{n-1} \ra W_n \ra P_{n-1}$, as
one sees using the octahedral axiom.
In particular, $W_n$ is in $\Proj_n$ for each $n$ and so
our Adams resolution provides us with choices of
cofibre sequences $W_{n} \ra X \ra X_{n}$ with 
$W_{n}$ in $\Proj_n$ and $X \ra X_{n}$ in $\I^n$.

If we apply the functor $\cS(-,Y)_*$ for some object $Y$ we get an
exact couple, which we display in unraveled form:
\begin{equation} \label{eq:uec} 
  \begin{array}{c}
  \xymatrix@C-2.4pc@R-0.5pc{ \relax
   \cS(X,Y)_* \ar[rd] && \cS(X_1,Y)_* \ar[ll] \ar[rd] && \cS(X_2,Y)_* \ar[ll] \ar[rd] && \cS(X_3,Y)_* \ar[ll] \\
   & \cS(P_0,Y)_* \ar[ru] \circar[ru] && \cS(P_1,Y)_* \ar[ru] \circar[ru] && \cS(P_2,Y)_* \ar[ru] \circar[ru] %
  }
  \end{array}
  \! \cdots .
\end{equation}
(We write $\cS(X,Y)_n$ for $\cS(\Sigma^n X,Y)$.)
This exact couple leads to a spectral sequence that we call
the \dfn{Adams spectral sequence} associated to the projective class.
The filtration on $\cS(X,Y)$ is the $\I$-adic filtration,
\ie that given by intersecting with the powers $\I^n$ of the ideal $\I$.

A dual construction using an injective class also leads to a
spectral sequence.  For example, one can obtain the original
Adams spectral sequence in this way by taking $X = S^{0}$,
the injectives to be retracts of products of (de)suspensions of
mod 2 \EM spectra,
and the maps in the ideal to be those which induce the zero
map in mod 2 singular cohomology.

For many results we will need to assume that our projective
class ``generates'', as described in the following definition.

\begin{defn}\label{de:generate}
A projective class is said to \dfn{generate} if for each
non-zero $X$ there is a projective $P$ such that $\cS(P,X) \neq 0$,
or, equivalently, if the ideal $\I$ contains no non-zero identity maps.
A third equivalent way to say this is that a map $Y \ra Z$ is an
isomorphism if and only if it is sent to an isomorphism by
the functor $\cS(P,-)$ for each projective $P$.
\end{defn}

One can show that a projective class $(\Proj ,\I )$ generates if
and only if one of its derived projective classes $(\Proj_{n},\I^n)$
generates.

The following result uses terminology from~\cite{bo:ccss}.\vspace*{-5pt}

\begin{prop}\label{pr:condconv}
  Let $(\Proj, \I)$ be a projective class such that 
  $\I$ is closed under countable coproducts.
\enlargethispage{10pt}
  Then the Adams spectral sequence abutting to $\cS(X,Y)$
  is conditionally convergent for all $X$ and $Y$
  if and only if the projective class generates.
\end{prop}

\begin{proof}
  Assume that the projectives generate and
  consider the unraveled exact couple pictured in~\eqref{eq:uec}.
  There is an exact sequence
  \[ 0 \lra \invlim \cS(X_k, Y)_* \lra \prod \cS(X_k, Y)_*
       \lra \prod \cS(X_k, Y)_* \lra \invlimone \cS(X_k, Y)_* \lra 0 \]
  in which the middle map is induced by the $(1 - \text{shift})$ map
  $\Wedge X_k \ra \Wedge X_k$.
  This map induces the identity map under $\cS(P,-)$ for each
  projective $P$, since the shift map $\Wedge X_k \ra \Wedge X_{k+1}$ is
  in $\I$ by assumption.
  Therefore $\Wedge X_k \ra \Wedge X_k$ is an isomorphism,
  since the projectives generate.
  This shows that
  \[ \invlim \cS(X_k, Y)_* = 0 = \invlimone \cS(X_k, Y)_* , \]
  and this is what it means for the associated spectral sequence 
  to be conditionally convergent.

  It isn't hard to see that if the projective class does not 
  generate, for example if the identity map $X \ra X$ is in $\I$
  for some non-zero $X$,
  then the spectral sequence abutting to $\cS(X,X)$ isn't
  conditionally convergent.
\end{proof}

When our projective class generates, we can characterize the objects
in $\Proj_{n}$ by the behaviour of the Adams spectral sequence.

\begin{prop}\label{pr:collapse}
If $X$ has length at most $n$, then
for each $Y$, the Adams spectral sequence abutting to $\cS(X,Y)$
collapses at $E_{n+1}$ and
has $E_{n+1} = E_{\infty}$ concentrated in the first $n$ rows.
If the projective class generates and $\I$ is closed under countable coproducts, 
then the converse holds,
and the spectral sequence converges strongly.
\end{prop}

We index our spectral sequence with the ``Adams indexing'',
so that the rows contain groups of the same homological degree
and the columns contain groups of the same total degree.

\begin{proof}
If $X$ is in $\Proj_n$, then one can easily
see that each $X_{s}$ appearing in an Adams resolution of $X$
is also in $\Proj_n$.
Therefore, the $n$-fold composites $X_{s} \ra \cdots \ra X_{s+n}$
are each zero.
But this implies that for each $r > n$ the differential $d_r$ is zero,
and so $E_{n+1} = E_{\infty}$.
Moreover, if $X$ is in $\Proj_n$, then $\I^n(X,Y) = 0$.
But $E_{\infty }$ is the associated graded of this filtration,
and so it must be zero except in the first $n$ rows.

Now we prove the converse.  
If the projective class generates and 
$\I$ is closed under countable coproducts, then by Proposition~\ref{pr:condconv}
the Adams spectral sequence is conditionally convergent.
Therefore, if it collapses, it converges strongly.
And if $E_{\infty }$ is concentrated in the first $n$ rows,
then the $n$th stage of the filtration on $\cS(X,Y)$ must be zero.
That is, $\I^{n}(X,Y) = 0$.
This is true for all $Y$, so $X$ is in $\Proj_{n}$.
\end{proof}

 From an Adams resolution, one can form the sequence
\[ 
  \xymatrixcolsep{0.9pc} 
  \xymatrixrowsep{1.0pc}
  \xymatrix{
   0 && X \ar[ll] 
   && P_0 \ar[ll] && P_1 \ar[ll] \circar[ll] && P_2 \ar[ll] \circar[ll] && 
                                 \cdots \period \circar[ll] \ar[ll] 
} \]
This is a projective resolution of $X$ with respect
to the projective class,
\ie each $P_s$ is projective and the sequence is $\Proj$-exact
at each spot.
Therefore the $E_2$-term of the Adams spectral sequence
consists of the derived functors of $\cS(-,Y)$ applied to $X$.
By the usual argument, these are independent of the choice of
resolution, and we denote them by $\Ext^k(-,Y)$, $k \geq 0$.
In fact, it is easy to see that from the $E_2$-term onwards the 
spectral sequence is independent of the choice of Adams resolution.

We note the following facts about the derived functors $\Ext^k(-,Y)$.
First, there is no reason to suspect that $\Ext^0(-,Y) = \cS(-,Y)$.
Indeed, the kernel of the natural map $\cS(X,Y) \ra \Ext^{0}(X,Y)$
is $\I(X,Y)$, and differentials in the Adams spectral sequence can 
prevent this map from being surjective.
Second, it is clear that if $X$ has projective dimension $n$, then
the groups $\Ext^k(X,Y)$ vanish for $k > n$.  To prove
the converse, we need to assume that the projective class generates.

\begin{prop}\label{pr:pdim}
If the projective class generates, then $X$ has projective dimension
at most $n$ if and only if $\Ext^k(X,-)$ vanishes for all $k > n$.
\end{prop}

\begin{proof}
Assume that $\Ext^k(X,Y)$ vanishes for each $Y$ and each $k > n$.
Consider an Adams resolution
\[
  \xymatrixcolsep{2.40em} 
  \xymatrixrowsep{3.25em}
  \begin{array}{c}
  \hspace*{.5em}
  \xymatrix@!R0{ \relax
   \makebox[2em][r]{$X = X_0$}
       \ar[rr] && X_1 \ar[ld] \circar[ld] \ar[rr] 
               && X_2 \ar[ld] \circar[ld] \\
   & P_0 \ar[lu] && P_1 \ar[lu] 
  }
  \end{array}
  \cdots
  \begin{array}{c}
  \xymatrix@!R0{ \relax
                  X_n \ar[rr]
               && X_{n+1} \ar[ld] \circar[ld] \ar[rr]
               && X_{n+2} \ar[ld] \circar[ld] \ar[rr]
               && X_{n+3} \ar[ld] \circar[ld] \\
   & P_n \ar[lu] && P_{n+1} \ar[lu] 
  && P_{n+2} \ar[lu] 
  }
  \end{array}
  \cdots
\]
of $X$.
Since the functors $\Ext^{k}(X,-)$ vanish for $k > n$, the sequence
$P_{n} \la P_{n+1} \la P_{n+2}$ is exact at $P_{n+1}$ after
applying the functor $\cS(-,Y)$ for any $Y$.
In particular, the map $P_{n+1} \ra X_{n+1}$ factors through $P_{n}$
giving a map $P_{n} \ra X_{n+1}$.
The composite $P_{n+1} \ra X_{n+1} \ra P_{n} \ra X_{n+1}$
is equal to the map $P_{n+1} \ra X_{n+1}$.
This shows that the composite $X_{n+1} \ra P_{n} \ra X_{n+1}$
is sent to the identity by the functor $\cS(P,-)$ for any projective $P$.
But since the projective class was assumed to generate, this implies
that the composite $X_{n+1} \ra P_{n} \ra X_{n+1}$ is an isomorphism.
Therefore $X_{n+1}$ is a retract of $P_{n}$, and so $X_{n}$ is as well.
Thus $X_{n}$ is projective, and
\[ 
  \xymatrixcolsep{0.8pc} 
  \xymatrixrowsep{1.0pc}
  \xymatrix{
   0 && X \ar[ll] 
   && P_0 \ar[ll] && P_1 \ar[ll] \circar[ll] && \cdots \ar[ll] \circar[ll] 
   && P_{n-1} \ar[ll] \circar[ll] && X_{n} \ar[ll] \circar[ll] && 0 \ar[ll]
} \]
displays that $X$ has projective dimension at most $n$.

The converse was noted above.
\end{proof}

So when our projective class generates,
saying that $X$ has projective dimension at most $n$ is equivalent to
saying that for each $Y$ the $E_{2}$-term of the Adams spectral sequence 
abutting to $\cS(X,Y)$ is concentrated in the first $n+1$ rows.

Here is another result that illustrates the importance of 
assuming that our projective class generates.

\begin{prop}\label{pr:bound}
If the projective class generates, then
$\pdim X + 1$ is an upper bound for
the length of $X$.
\end{prop}

This follows from Proposition~\ref{pr:collapse} in the case that $\I$
is closed under countable coproducts.

\begin{proof}
If
\[ 
  \xymatrixcolsep{0.9pc} 
  \xymatrixrowsep{1.0pc}
  \xymatrix{
   0 && X \ar[ll]
  && P_0 \ar[ll]             && P_1    \ar[ll] \circar[ll] 
  && P_2 \ar[ll] \circar[ll] && \cdots \ar[ll] \circar[ll]
} \]
is a projective resolution of $X$ with respect to the projective
class then it can be filled out to an Adams resolution
\[
  \xymatrixcolsep{0.6pc} 
  \xymatrixrowsep{1.05pc}
  \begin{array}{c}
  \xymatrix{ \relax
   \makebox[2em][r]{$X = X_0$}
   \ar[rr] && X_1 \ar[ld] \circar[ld] \ar[rr] && X_2 \ar[ld] \circar[ld] \ar[rr] && X_{3} \ar[ld] \circar[ld] \\
   & P_0 \ar[lu] && P_1 \ar[lu] && P_2 \ar[lu]
} 
  \end{array}
  \cdots
\]
in which the composites $P_n \ra X_n \ra P_{n-1}$ 
equal the maps $P_n \ra P_{n-1}$ appearing in the first diagram.
If the projective resolution is finite, say with $P_n = 0$ for $n > k$,
then $\cS(P,X_n) = 0$ for all projectives $P$ and all $n > k$.
Because the projective class generates, $X_n = 0$ for $n > k$.
Thus $X_k$ has length at most one, $X_{k-1}$
has length at most two, and inductively, $X = X_0$ has length at
most $k+1$, completing the argument.
\end{proof}

In general the projective dimension of $X$ will be larger than
its length;  the difference comes about because when we measure the
length by building up $X$ using projectives, we don't insist that the 
connecting maps $X_k \ra X_{k+1}$ be in $\I$.

\mysection{Abstract phantom maps}\label{se:phantoms}

In this section we discuss phantom maps in an axiomatic setting.  
We begin in the first part by defining phantom maps and describing some
assumptions that we will need to state our results.
The second part is a short study of various flavours of weak colimits.
This is essential material for the third part of the section,
which gives our results on phantom maps.
This section is based on joint work with Neil Strickland~\cite{chst:pmht}.

We remind the reader that $\cS$ denotes a triangulated category
having all set-indexed coproducts.

\subsection{The phantom projective class} \label{se:ppc}

We begin with a definition.

\begin{defn}\label{de:finite}
  An object $W$ in $\cS$ is \dfn{finite} if for any set-indexed
  collection $\{X_\alpha\}$ of objects of $\cS$, the natural map
  \[ \bigoplus_\alpha \cS(W,X_\alpha) \lra \cS(W,\bigWedge_\alpha X_\alpha) \]
  is an isomorphism.
\end{defn}

To illustrate that this is a reasonable definition, we describe
the finite objects in
the categories that we study in the last few sections of the paper.
In the stable homotopy category, an object is finite if and only if
it is isomorphic to a (possibly desuspended) suspension spectrum of
a finite CW-complex.
In the derived category of a ring, an object is finite if and only if
it is isomorphic to a
bounded complex of finitely generated projectives.
In both cases, the finite objects are precisely those that can
be built from a finite number of copies of the spheres
($S^{n}$ and $\Sigma^{n} R$, respectively) using cofibres
and retracts.

Write $\Proj $ for the collection of retracts of wedges of finite objects.

A map $X \ra Y$ is a \dfn{phantom map} if for each finite $W$ and each
map $W \ra X$ the composite $W \ra X \ra Y$ is zero.
Write $\I $ for the collection of phantom maps.

\begin{defn}\label{de:gppc}
  We say that $\cS$ \dfn{has a phantom projective class} if
  $(\Proj,\I)$ is a projective class.
  We say that $\cS$ \dfn{has a generating phantom projective class} if
  it has a phantom projective class and this projective class generates.
\end{defn}

That the projective class generates
says that if $\cS(W,X) = 0$ for each finite $W$, then $X = 0$.
In other words, this says that $\cS $ is \dfn{compactly generated}, 
in the terminology of Neeman~\cite{ne:tba}.

Assuming that $(\Proj,\I)$ is a projective class
is equivalent to assuming that for each
$X$ there is a set $\{X_{\alpha }\}$ of finite objects such
that every map $W \ra X$ from a finite object to $X$ factors
through some $X_{\alpha }$.
In particular, if there is a set of isomorphism classes of
finite objects, then $\cS$ has a phantom projective class.
Thus we have replaced a set-theoretic condition
with the slightly more general and natural condition
that $(\Proj,\I)$ be a projective class.

The stable homotopy category and the derived category of
a ring are examples of triangulated categories with generating
phantom projective classes.

Our strongest results will be possible when $\cS$ is a
``Brown'' category:

\begin{defn}\label{de:br}
We say that $\cS$ is a \dfn{Brown} category if the
following holds for any two objects $X$ and $Y$.
Regarding $\cS(-,X)$ and $\cS(-,Y)$ as functors from finite objects to
abelian groups, any natural transformation $\cS(-,X) \ra \cS(-,Y)$
is induced by a map $X \ra Y$.
\end{defn}

In familiar settings, this can be rephrased as the assumption that
natural transformations between representable homology theories are
representable.
The stable homotopy category and 
the derived category of a countable ring
are Brown categories. 
(See \cite[Theorem~4.1.5]{hopast:ash} or \cite[Section~5]{ne:tba}.)

In order to prove results about the phantom projective class,
we need a digression on weak colimits.

\subsection{Weak colimits} \label{se:wc}

Colimits rarely exist in a triangulated category, so in this section
we introduce weaker variants that turn out to be quite useful. 

\begin{defn}
  A category $\cC$ is \dfn{small} if its class of objects forms a set.
  A \dfn{diagram} in $\cS$ is a functor $F$ from a small category $\cC$
  to $\cS$.
  A \dfn{cone} from a diagram $F$ to an object $X$ is a natural
  transformation from $F$ to the constant diagram at $X$.
  In other words, for each $\alpha$ in $\cC$ we are given a
  map $i_{\alpha }: F(\alpha ) \ra X$ such that for each
  map $\alpha \ra \beta $ in $\cC$ the triangle
  \[ \xymatrixcolsep{1.3pc} 
     \xymatrixrowsep{0.25pc}
     \xymatrix{ \relax
     F(\alpha ) \ar[rd] \ar[dd]  \\
     & X \\
     F(\beta ) \ar[ru]
  } \]
  commutes.
  A \dfn{weak colimit}
  of a diagram $F$ is a cone through which every other cone factors.
  If we require the factorization to be unique, this is
  the definition of a \dfn{colimit}.
\end{defn}

If $F$ is a diagram, then a weak colimit of $F$ always exists.
It will not be unique, but there is a distinguished
choice of weak colimit defined in the following way.
Writing $X_{\alpha }$ for $F(\alpha )$, there is a natural map
\[ \bigWedge_{\alpha \ra \beta} X_\alpha \lra \bigWedge_{\gamma} X_\gamma . \]
The first coproduct is over the non-identity morphisms of $\cC$, and the
second is over the objects.  
(We omit the identity morphisms for reasons explained in Note~\ref{no:non-id}.)
The restriction of this map to the summand $X_\alpha$ indexed 
by $\alpha \ra \beta$ is
\[ X_\alpha \xra{(1,-F(\alpha \ra \beta))} X_\alpha \Wedge X_\beta
   \xra{\text{inclusion}} \bigWedge_\gamma X_\gamma . \]
Let $X$ be the cofibre of the natural map, so that $X$ sits
in a cofibre sequence
\[ \bigWedge_{\alpha \ra \beta} X_\alpha \lra \bigWedge_{\gamma} X_\gamma 
    \lra X \lra \bigWedge_{\alpha \ra \beta} \Sigma X_\alpha . \]
The map $\bigWedge_\gamma X_\gamma \ra X$ gives a cone $i$ to $X$,
and it is easily checked that $X$ and $i$ form a weak colimit of the 
diagram $F$.
Also, given another weak colimit $(X', i')$ constructed in the same way, there
is an isomorphism $h : X \ra X'$ such that $ih=i'$.
(This isomorphism might not be unique.)
We call $X$ the \dfn{standard weak colimit} of $F$.

The virtue of standard weak colimits is that they always exist
and are easy to describe.  
However, in many cases they are too large, and there are
other types of cones that are more useful.

\begin{defn}\label{de:minimal}
A cone $i : F \ra X$ is a \dfn{minimal cone} if the natural map
\[ \dirlim \cS(W,X_{\alpha }) \ra \cS(W,X) \]
is an isomorphism for each finite $W$.
We say that a cone to $X$ is a \dfn{minimal weak colimit} 
if it is a minimal cone and a weak colimit.
\end{defn}

When they exist, minimal weak colimits behave well, as we will see
in the following proposition.
The reason for introducing the weaker notion of a minimal cone
is that one is often able to verify that a cone is minimal without
being able to prove that it is a weak colimit.  And we will find
that minimal cones share some of the nice properties of minimal
weak colimits.

\begin{prop}\label{pr:cones}
\begin{roenumerate}
\item Assume that $\cS$ has a phantom projective class.
Then any object $X$ is a minimal cone on a filtered diagram of finite
objects.
\item Assume that $\cS$ is a Brown category.
Then any minimal cone on a diagram of projective objects is a weak colimit.
\item Assume that the finite objects generate, \ie that for each
non-zero $X$ there is a finite object $W$ and a non-zero map $W \ra X$.
Then a minimal weak colimit is unique up to (non-unique) isomorphism and 
is a retract of any other weak colimit.
\end{roenumerate}
\end{prop}

\begin{proof}
We begin by proving (\emph{i}).
Choose a set $\{X_{\alpha}\}$ of finite objects such
that every map $W \ra X$ from a finite object to $X$ factors
through some $X_{\alpha}$.
Consider the thick subcategory $\cC$ generated by the $X_{\alpha}$.
In other words, consider the smallest full subcategory that contains
each $X_{\alpha}$ and is closed under taking cofibres, desuspensions and
retracts.
Since the collection of $X_{\alpha}$'s is a set, one can show
that there is a set $\cC'$ containing a representative of each
isomorphism class of objects in $\cC$.
Let $\Lambda(X)$ denote the category whose objects are
the maps $W \ra X$ with $W$ in $\cC'$ and whose morphisms
are the obvious commutative triangles.
One can check that this is a filtered category, using the fact that
$\cC$ is a thick subcategory.
There is a natural functor $\Lambda(X) \ra \cS$ sending an
object $W \ra X$ to $W$, and there is a natural cone from
this diagram to $X$.
Consider a finite object $V$.  We must show that the
natural map
\[
\underset{W \ra X \in \Lambda(X)}{\dirlim} \cS(V,W) \ra \cS(V,X)
\]
is an isomorphism.
Surjectivity is easy, since any map from $V$ to $X$ factors through
some $X_{\alpha }$.
Now we prove injectivity.  Since the diagram is filtered,
a general element of the colimit can be represented by a map $V \ra W$
for some $W \ra X \in \Lambda(X)$.
Suppose that such an element is sent to zero, \ie that the
composite $V \ra W \ra X$ is zero.  Define $V'$ to be
the cofibre of the map $V \ra W$.
The map from $W$ to $X$ factors through $V'$.
Now $V'$ might not be in $\cC$, but the map from it to $X$ factors through
an object $W'$ in $\cC$.
The composite $V \ra W \ra V' \ra W'$ is zero, since the first
three terms form a cofibre sequence, so the element of $\cS(V,W)$
given by the map $V \ra W$ goes to zero under the map $W \ra W'$.
Thus $V \ra W$ represents the zero element of the colimit,
and we have proved injectivity.

Next we prove (\emph{ii}).
Let $X$ be a minimal cone on a diagram $\{X_{\alpha }\}$ of projective objects.
Given another cone $Y$ we must show that
there exists a map $X \ra Y$ commuting with the cone
maps.
To construct this map, we will
produce a natural transformation from $\cS(-,X)$ to
$\cS(-,Y)$, regarded as functors on finite objects.
Since there is a cone to $Y$, there is a natural transformation
$\dirlim \cS(-,X_{\alpha }) \ra \cS(-,Y)$.
And since $X$ is a minimal cone, we have a natural isomorphism
$\dirlim \cS(-,X_{\alpha }) \ra \cS(-,X)$.
The composite of the inverse of the isomorphism with the
map to $\cS(-,Y)$ is the natural transformation we said
that we would produce.
Thus there is a map $X \ra Y$ inducing this natural
transformation, as we have assumed that $\cS$ is a Brown category.
This is the map we seek.  By construction, the triangles
commute up to phantom maps.  But since we assumed that each $X_{\alpha}$
is projective, the triangles actually commute.

Finally, we prove (\emph{iii}).
Let $\{X_{\alpha }\}$ be a diagram with minimal weak colimit $X$ 
and minimal cone $Y$.
Because $X$ is a weak colimit, there is a map $X \ra Y$ commuting
with the cone maps.
Because both cones are minimal, this map is an isomorphism under
$\cS(W,-)$ for each finite $W$. 
Since the finite objects generate, we can conclude that 
the map $X \ra Y$ is an isomorphism.
In particular, if $X$ and $Y$ are both minimal weak colimits,
then they are isomorphic.

For the second part of (\emph{iii}), assume that $X$ is a minimal
weak colimit and that $Y$ is a weak colimit of a diagram
$\{X_{\alpha }\}$.
Since both cones are weak colimits, there are maps
$X \ra Y$ and $Y \ra X$ commuting with the cone maps.
Because $X$ is minimal, the composite $X \ra Y \ra X$ is
an isomorphism, again using that the finite objects generate.
\end{proof}

Minimal weak colimits earned their name by way of part (\emph{iii}) of
this Proposition.
As far as I know, they were first defined by 
Margolis~\cite[Section~3.1]{ma:ssa}.

Next we consider some examples.  Let $\cC$ be the non-negative integers,
where there is one map $m \ra n$ when $m \leq n$, and no maps otherwise.
A functor $F : \cC \ra \cS$ is a diagram of the form
\[ X_0 \lra X_1 \lra X_2 \lra \cdots . \]
The minimal weak colimit is the cofibre of the usual $(1-\text{shift})$ map 
\[ \bigWedge X_k \lra \bigWedge X_k . \]
Thus the minimal weak colimit is what is usually called the
\dfn{telescope} of the sequence.
The standard weak colimit is in this case much less manageable.

Our second example concerns the weak pushout.  

\begin{lemma} \label{le:wpo}
  Given a commutative diagram
  \[ \xymatrix{ \relax
       & V \ar@{=}[r] \ar[d] & V \ar[d] \\
H \ar@{=}[d] \ar[r] & W \ar[r] \ar[d] & X \ar[r] \ar[d] & \Sigma H \ar@{=}[d]\\
H            \ar[r] & Y \ar[r] \ar[d] & Z \ar[r] \ar[d] & \Sigma H           \\
       & \Sigma V \ar@{=}[r] & \Sigma V 
  } \]
  with exact rows and columns, the centre square is both a weak pushout
  and a weak pullback.              
  Moreover, the standard weak pushout of the diagram
  \[ \xymatrix{ \relax
     W \ar[r] \ar[d] & X \\
     Y 
  } \]
  fits into a diagram of the above form.
\end{lemma}

The proof is omitted, as we will not make use of this result.

\begin{note}\label{no:non-id}
  In the definition of the standard weak colimit, the coproduct
\[
\bigWedge_{\alpha \ra \beta} X_\alpha
\]
  was taken over the \emph{non-identity} morphisms of the indexing
  category.  Had we instead taken the coproduct over all of the
  morphisms, we would have obtained a different distinguished weak colimit, 
  but the second half of Lemma~\ref{le:wpo} would no longer be true.
\end{note}

\begin{note}
  For a weak pushout it turns out that the minimal weak colimit is
  less useful than the standard weak colimit.  For example, the
  standard weak colimit of the diagram
  \[ \xymatrix{ \relax
     W \ar[r] \ar[d] & 0 \\
     0 
  } \]
  is $\Sigma W$, while the minimal weak colimit is $0$.
\end{note}

\subsection{Consequences}\label{se:cons}

With the work we have done, we can immediately prove the
following theorem.

\begin{thm}\label{th:limi}
Assume that $\cS$ has a phantom projective class.
If $X$ is a minimal cone on a filtered diagram $\{X_{\alpha}\}$ 
of finite objects, then the $E_{2}$-term of the phantom Adams spectral
sequence abutting to $\cS(X,Y)$ is given by
\[
E_{2}^{s} = \invlims \cS(X_{\alpha },Y) ,
\]
where $\invlims$ denotes the $s$th derived functor of 
the inverse limit functor.
Here $s$ is the homological degree and we have suppressed the internal degree.
\end{thm}

\begin{proof}
  We prove this by constructing a specific Adams resolution of $X$ with
  respect to the phantom projective class.
  Consider the sequence
  \[ 0 \lla X \lla \bigWedge_\alpha X_\alpha \lla \bigWedge_{\alpha \ra \beta} X_\alpha \lla \bigWedge_{\alpha \ra \beta \ra \gamma} X_\alpha \lla \cdots . \]
  The wedges are over sequences of morphisms in the category over which
  the diagram $\{X_\alpha\}$ is indexed (and here identity morphisms
  are included).
  Because $X$ is a weak colimit of the $X_\alpha$, there are given
  maps $i_\alpha : X_\alpha \ra X$.
  The map $\Wedge_{\alpha} X_\alpha \ra X$ 
  is equal to $i_\alpha$ on the $\alpha$ summand.
  The map $\Wedge_{\alpha \ra \beta} X_\alpha \ra \Wedge_{\alpha} X_\alpha$
  sends the summand $X_\alpha$ indexed by the map $\alpha \ra \beta$
  to the $\alpha$ summand of the target using the identity map and
  to the $\beta$ summand of the target using the negative of the map 
  $X_{\alpha} \ra X_{\beta}$.
  In general, one gets an alternating sum.
  When we apply the functor $\cS(W,-)$ for finite $W$
  we get the sequence used for computing the derived functors of 
  $\dirlim$~\cite[Appendix II, Section 3]{gazi:cfht}.
  Since these vanish (because we have a filtered colimit) and since
  $\dirlim \cS(W,X_{\alpha}) = \cS(W,X)$ (because we have a minimal cone), 
  the sequence obtained is exact.  That is, the sequence above is a phantom
  projective resolution of $X$.  Therefore it is part of an Adams tower
  for $X$.
  The $E_2$-term of the Adams spectral sequence
  obtained by applying $\cS(-,Y)$ is the cohomology of the sequence
  \[ 0 \lra \bigoplus_\alpha \cS(X_\alpha,Y) \lra \bigoplus_{\alpha \ra \beta} \cS(X_\alpha,Y) \lra \cdots . \]
  But the cohomology at the $s$th place is just $\invlims \cS(X_\alpha, Y)$,
  again by~\cite[App.\ II]{gazi:cfht}.  
\end{proof}

Under the assumptions described in the first part of the section, 
the phantom projective class is very well behaved.

\begin{thm}\label{th:pdim1}
If $\cS$ is a Brown category with a generating phantom projective
class, then any object $X$ has projective dimension
at most one.
\end{thm}

This result was proved independently by several 
people~\cite{oh:vtaybkss,ne:tba,chst:pmht}.
The proofs by Neeman and Christensen-Strickland are
essentially the same, and both are phrased in an axiomatic
setting similar to that presented here.
On the other hand, Ohkawa's proof was written in the context
of the stable homotopy category and makes use of CW-structures,
so it is not clear whether it goes through in the same generality.
Ohkawa also noticed the consequence for the Adams spectral
sequence, which is our Corollary~\ref{co:milnor}.

\begin{proof}
Let $X$ be an object of $\cS$.
We saw in Proposition~\ref{pr:cones}~(\emph{i}) that $X$ is a minimal cone 
on a filtered diagram $\{X_{\alpha} \}$ of finite objects.  
The standard weak colimit $Y$ of this diagram
lies in a cofibre sequence
\[ \bigWedge_{\alpha \ra \beta} X_\alpha \lra \bigWedge_{\gamma} X_\gamma 
    \lra Y \lra \bigWedge_{\alpha \ra \beta} \Sigma X_\alpha . \]
By Proposition~\ref{pr:cones}~(\emph{ii}) and~(\emph{iii}), $X$ is a retract of $Y$.
This implies that the fibre $P$ of the natural map
$\Wedge_{\gamma } X_{\gamma } \ra X$ is a retract of
$\Wedge_{\alpha  \ra \beta } X_{\alpha }$ and thus is projective.
Moreover, the connecting map $X \ra \Sigma P$ is phantom because $X$ 
is a minimal cone, so we have constructed a projective resolution
\[
0 \lra P \lra Q \lra X \lra 0 .
\]
Therefore, $X$ has projective dimension at most one.
\end{proof}
 
\begin{cor}\label{co:milnor}
If $\cS$ is a Brown category with a generating phantom projective
class, then the Adams spectral sequence collapses at the $E_{2}$-term and the
composite of two phantom maps is zero.
Moreover, if $\{X_\alpha\}$ is a filtered diagram of finite objects with
minimal cone $X$,
then there is a short exact sequence
\[ 0 \lra \invlimone \cS(\Sigma X_\alpha, Y) \lra \cS(X,Y) 
     \lra \invlim \cS(X_\alpha, Y) \lra 0 \]
natural in $Y$.
The kernel consists of the phantom maps from $X$ to $Y$,
and $\invlimi \cS(\Sigma X_\alpha, Y)$ is zero for $i \geq 2$.
\end{cor}

\begin{proof}
Consider the Adams spectral sequence abutting to $\cS(X,Y)$.
By Theorem~\ref{th:pdim1}, $X$ has projective dimension at most one.
Since the $E_2$-term of the Adams spectral sequence consists of the
derived functors of $\cS(-,Y)$ applied to $X$, it must vanish in
all but the first two rows.
There being no room for differentials, it collapses at $E_{2}$,
degenerating into the displayed Milnor sequence.

By Proposition~\ref{pr:bound}, we see that each object has length at most two,
and so the composite of two phantoms must be zero.
\end{proof}

\mysection{Topological phantom maps}\label{se:top-phantoms}

In this section $\cS$ denotes the stable homotopy category and
we usually write $[X,Y]$ for the set of morphisms from
$X$ to $Y$.
There are many descriptions of this category;
a good one can be found in the book by Adams~\cite{ad:shgh}.

There are three parts to this section.  In the first we discuss
phantom maps, \ie maps which are zero when restricted to any finite
spectrum, and we use the results of the previous section
to conclude that there is a generalized Milnor sequence.
In the second part of this section 
we discuss skeletal phantom maps, \ie maps which
are zero when restricted to each skeleton of the source.
And in the third we discuss superphantom maps, \ie maps which are
zero when restricted to any (possibly desuspended) suspension
spectrum.

\subsection{Phantom maps and a generalized Milnor sequence}\label{se:gms}

In the stable homotopy category, a finite spectrum is
one isomorphic to a (possibly desuspended) suspension spectrum
of a finite CW-complex. %
As in the previous section, a map $f : X \ra Y$ is said to be \dfn{phantom}
if the composite $W \ra X \ra Y$ is zero for each finite spectrum
$W$ and each map $W \ra X$.
Phantom maps form an ideal which we denote $\I$.
Write $\Proj$ for the collection of all retracts of wedges
of finite spectra.
One can show that there exists a set $\cF'$ of finite spectra containing a 
representative of each isomorphism class~\cite[Prop.~3.2.11]{ma:ssa}, and
it follows that $(\Proj,\I)$ is a projective class.
Since $\Proj$ contains the spheres, it is in fact a generating
projective class.
Also, it is well-known that the stable homotopy category
is a Brown category (Definition~\ref{de:br}).

We should point out that there do exist non-zero phantom maps.
For example, if $G$ is any non-zero divisible abelian group,
then the Moore spectrum $S(G)$ is not a retract of a wedge 
of finite spectra.  If it was, then applying integral
homology would show that $G$ is a retract of a sum of finitely
generated abelian groups, which is impossible.  So there is
a non-zero phantom map with source $S(G)$.

For another example, consider the natural map 
\[
\bigWedge X_{\alpha} \lra \prod X_{\alpha }
\]
for some indexed collection of spectra.
For finite $W$, $[W,\Wedge X_{\alpha }] = \oplus [W,X_{\alpha }]$,
so we get a monomorphism
\[
[W,\bigWedge X_{\alpha}] \lra [W,\prod X_{\alpha }] .
\]
Thus the fibre of the map $\Wedge X_{\alpha } \ra \prod X_{\alpha }$
is phantom.
It is non-zero if and only if the map from the wedge to the product
is not split.
This is the case for the map
\[
\bigWedge H\Z \lra \prod H\Z 
\]
from the countable wedge of integral \EM spectra
to the countable product.
Indeed, if this map splits, then so does the map
\[
\bigoplus \Z \lra \prod \Z 
\]
of abelian groups.  But the cokernel of this map contains the
element $[1,2,4,8,\dots ]$.  This element is divisible by all
powers of $2$, so the cokernel is not a subgroup of the product.

Other examples will appear in the second and third parts of this section.

The stable homotopy category is a Brown category with a
generating phantom projective class.
Thus, the results of Section~\ref{se:cons} hold and we find
that every spectrum has projective dimension at most one
and length at most two, and that the
composite of two phantom maps is zero.
Also, a minimal cone on a diagram of projective objects 
is automatically a minimal weak colimit
and leads to a generalized Milnor sequence.

Minimal cones arise in practice.
If a diagram is filtered, then to check
that a cone is minimal, it suffices to check that it becomes
a colimiting cone in homotopy groups.
This is proved using the fact that a finite spectrum is
built from a finite number of spheres using cofibres.
One uses induction on the number of cells, that filtered
colimits are exact, and the five-lemma.
For this reason, all of our examples will involve filtered diagrams.

First of all, if $X$ is a CW-spectrum (in the sense of Adams~\cite{ad:shgh})
and $\{X_{\alpha } \}$ is a 
filtered collection of finite CW-subspectra whose union is $X$, then
$X$ is the minimal weak colimit of the $X_{\alpha }$.

For another example, let $\{G_{\alpha } \}$ be a filtered 
diagram of abelian groups with colimit $G$.  Then the 
\EM spectrum $H G$ is a minimal cone on the 
diagram $\{H G_{\alpha } \}$.
(By construction, this cone becomes a colimiting cone in homotopy groups.)
I don't believe that such a minimal cone is always a minimal weak colimit.

More generally, a filtered homotopy colimit of spectra
(taken in some geometric category of spectra) is a minimal cone,
because homotopy groups commute with filtered homotopy colimits.
If the spectra in the diagram are projective, then the homotopy
colimit is a weak colimit.

 From the first example and Corollary~\ref{co:milnor}, 
we can state the following theorem.

\begin{thm}\label{th:milnor}
  Let $X$ be a CW-spectrum and let $\{X_\alpha\}$ be a filtered diagram
  of finite CW-subspectra whose union is $X$.
  For any spectrum $Y$ there is a short exact sequence
  \[ 0 \lra \invlimone [\Sigma X_\alpha, Y] \lra [X,Y] 
       \lra \invlim [X_\alpha, Y] \lra 0 . \]
  The kernel consists precisely of the phantom maps.
  Moreover, $\invlimi [\Sigma X_\alpha, Y]$ vanishes for $i \geq 2$.   \qed
\end{thm}

This generalizes results of Pezennec~\cite{pe:ptff},
Huber and Meier~\cite{hume:ctic}, and Yosimura~\cite{yo:hcbpc}.
Pezennec makes the assumption that $Y$ has finite type.
Huber and Meier make the weaker assumption that the cohomology
theory represented by $Y$ is related by a universal coefficient
sequence to a homology theory of finite type, while
Yosimura drops the assumption that the homology theory has finite type.
Our point is that no restriction on $Y$ is necessary; this was proved
independently, and earlier, by Ohkawa~\cite{oh:vtaybkss}.

Using a slightly more elaborate proof (and a slightly different
projective class), one can show that the assumption
that each $X_\alpha$ is finite can be replaced by the assumption that
there are no phantom maps from $X_\alpha$ to $Y$ for each $\alpha$.

\subsection{Skeletal phantom maps}

There is a related but smaller ideal of maps which have also been
called phantom maps in the literature.  We begin with some background
on cellular towers.
We treat cellular towers in this abstract manner because we want to
define and use them without stepping outside of the homotopy category.
One reason for this is that we want to make it clear that our
results do not depend on a particular choice of model for the
category of spectra.
But more importantly, we would like our arguments to go through
in any nice enough triangulated category.

\begin{defn} \label{de:cellular-tower}
  Let $X$ be a spectrum.  A \dfn{cellular tower} for $X$ is a diagram
  \begin{equation} \label{eq:cellular} 
  \begin{array}{c} \xymatrix{ \relax
  \cdots \ar[r] & X^{(n)} \ar[r] \ar[drrr] & X^{(n+1)} \ar[r] \ar[drr] &
                  X^{(n+2)} \ar[r] \ar[dr] & \cdots \\
  &&&& X
  } 
  \end{array}
  \end{equation}
  satisfying:
  \begin{roenumerate}
    \item $X$ is the telescope of the sequence 
          $\cdots \ra X^{(n)} \ra X^{(n+1)} \ra \cdots$.
    \item The fibre of the map $X^{(n)} \ra X^{(n+1)}$ is a 
          wedge of $n$-spheres.
    \item The inverse limit of abelian groups $\invlim H_*(X^{(n)})$
          is zero, where $H_*$ denotes integral homology.
  \end{roenumerate}
  We say that $X^{(n)}$ is an \mdfn{$n$-skeleton} of $X$.
\end{defn}

The first condition says that the sequence
\[ \Wedge X^{(n)} \lra \Wedge X^{(n)} \lra X \]
is a cofibre sequence, where the first map is the $(1 - \text{shift})$ map.

The above definition is taken from Margolis' book~\cite[Section~6.3]{ma:ssa}, 
which is also the source of the results below whose proofs are omitted.

\begin{prop}
  A diagram~\eqref{eq:cellular} is a cellular tower for $X$ if and
  only if all of the following conditions hold:
  \begin{roenumerate}
  \item The map $\pi_i(X^{(n)}) \ra \pi_i(X)$ is an isomorphism for
    each $i < n$.
  \item Each $H_n(X^{(n)})$ is a free abelian group and 
    $H_n(X^{(n)}) \ra H_n(X^{(n+1)})$ is an epimorphism.
  \item $H_i(X^{(n)}) = 0$ for $i > n$.                               \qed
  \end{roenumerate}      
\end{prop}

\begin{prop}
  Each spectrum $X$ has a cellular tower.         \qed
\end{prop}

\begin{prop} \label{pr:nat}
  Let $\cdots \ra X^{(n)} \ra X^{(n+1)} \ra \cdots \ra X$ be a cellular
  tower for $X$ and 
  let $\cdots \ra Y^{(n)} \ra Y^{(n+1)} \ra \cdots \ra Y$ be a cellular
  tower for $Y$.
  Given any map $X \ra Y$, there exist maps $X^{(n)} \ra Y^{(n)}$ making the
  following diagram commute:
  \[ \xymatrix{ \relax
  \cdots \ar[r] & X^{(n)} \ar[r] \ar[d]
    & X^{(n+1)} \ar[r] \ar[d] & \cdots \ar[r] & X \ar[d] \\
  \cdots \ar[r] & Y^{(n)} \ar[r] & Y^{(n+1)} \ar[r] & \cdots \ar[r] & Y 
    \period
  } \]
  \\*[-.25in]
  \qed
\end{prop}

\begin{defn}
  We say that a spectrum $X$ is an \mdfn{$n$-skeleton} if it is an $n$-skeleton
  of some spectrum $Y$.
  We say that $X$ is a \dfn{skeleton}, or is \dfn{skeletal}, if $X$ is
  an $n$-skeleton for some $n$.
\end{defn}

It is easy to see that if $X$ is an $n$-skeleton, then $X$ is an
$n$-skeleton of itself.

\begin{prop}
  A spectrum $X$ is an $n$-skeleton if and only if $H_i(X) = 0$ for
  $i > n$ and $H_n(X)$ is a free abelian group (possibly zero).
  In particular, $X$ is skeletal if and only if it has bounded
  above integral homology.     \qed 
\end{prop}

We call a map $f : X \ra Y$ a \dfn{skeletal phantom map} if the composite
$W \ra X \ra Y$ is zero for each skeleton $W$
and each map $W \ra X$. 
By Proposition~\ref{pr:nat}, it suffices to test this for the skeleta
$X^{(n)}$ of a fixed cellular tower for $X$.
This shows that skeletal phantoms form part of a projective class,
because it allows us to restrict to a \emph{set} of test objects.
In this case the projectives are retracts of wedges of skeletal spectra.
In the cofibre sequence $\Wedge X^{(n)} \ra X \ra \Wedge \Sigma X^{(n)}$,
$\Wedge X^{(n)}$ is projective and $X \ra \Wedge \Sigma X^{(n)}$ is
a skeletal phantom.

Since skeletal phantoms are phantoms, it follows from the previous
part of this section that the composite of two skeletal phantoms is zero.
In fact, this is obvious because the cofibre sequence
$\Wedge X^{(n)} \ra X \ra \Wedge \Sigma X^{(n)}$ shows that every
$X$ has length at most two.
That the composite of two skeletal phantoms is zero has been
known for some time.  (See~\cite{gr:oph} and~\cite{grmc:upm}.)

Applying the functor $[-,Y]$ to the cofibre sequence 
$\Wedge X^{(n)} \ra X \ra \Wedge \Sigma X^{(n)}$
immediately reveals another Milnor sequence
\[ 0 \lra \invlimone [\Sigma X^{(n)}, Y] \lra [X,Y] 
       \lra \invlim [X^{(n)}, Y] \lra 0 . \]
It is clear that the kernel consists precisely of the skeletal phantoms,
and that, for $i \geq 2$, $\invlimi [X^{(n)},Y] = 0$ (since, for $i \geq 2$,
$\invlimi$ is zero for a diagram indexed by the integers).

If $X$ is a spectrum with finite skeleta, then a map $X \ra Y$ is
phantom if and only if it is a skeletal phantom.  But in general the
ideal of skeletal phantoms is strictly smaller than the ideal of
phantom maps.  For example, we saw above that for $G$ a non-zero
divisible abelian group, 
the Moore spectrum $S(G)$ is the source of a non-zero phantom map.
But $S(G)$ is skeletal, since it has bounded above integral homology,
and so $S(G)$ is not the source of a non-zero skeletal phantom.

On the other hand, there are non-zero skeletal phantoms.
Consider $\Hp$, the mod $p$ \EM spectrum.
This spectrum has finite skeleta, so it suffices to show that
there is a non-zero phantom with source $\Hp$.
This is equivalent to showing that $\Hp$ is not a retract
of a wedge of finite spectra.
If it was and $\Hp \ra \Wedge W_\alpha$ was the inclusion,
then we would look at the composite $\Hp \ra \Wedge W_\alpha
\ra \prod W_\alpha$, which is a monomorphism in homotopy groups.
But any such map is zero because there are no maps from
$\Hp$ to a finite spectrum (see~\cite{ma:ems}, \cite{li:dems},
or~\cite{ra:lrc}) and so we would conclude that $\Hp$ has no
homotopy.  
This type of argument will appear repeatedly in what follows.

A non-zero phantom from $\Hp$ is also an example of a phantom
map which is not divisible by $p$, since $p$ kills $[\Hp,Y]$
for any $Y$.

\subsection{Superphantom maps}

There is another special class of phantom maps.
We call a map $f : X \ra Y$ a \dfn{superphantom map} if the
composite $W \ra X \ra Y$ is zero for each (possibly desuspended)
suspension spectrum $W$ and each map $W \ra X$.
Again the superphantoms form an ideal which is part of a 
projective class.  The projectives here are all retracts of
wedges of (possibly desuspended) suspension spectra.
To see that we have a projective class, one uses the
following lemma, which allows us to use a \emph{set} of objects
in order to test whether a map $X \ra Y$ is superphantom.

\begin{lemma}
  A map $X \ra Y$ is a superphantom if and only if for each $n$ the composite
  \[ \Sigma^{-n} \Sigma^{\infty} \Omega^{\infty} \Sigma^{n} X \lra X \lra Y\]
  is zero.
\end{lemma}

Recall that $\Sigma^\infty$ is left adjoint to $\Omega^\infty$.
The map $\Sigma^{-n} \Sigma^{\infty} \Omega^{\infty} \Sigma^{n} X \ra X$
is the $n$th desuspension of the counit of the adjunction.

\begin{proof}
  The ``only if'' direction is clear.  So suppose $f : X \ra Y$ is
  a map such that each composite
  $\Sigma^{-n} \Sigma^{\infty} \Omega^{\infty} \Sigma^{n} X \ra X \ra Y$
  is zero, and let $W$ be a space.
  Consider a map $\Sigma^{-n} \Sigma^\infty W \ra X$.
  This map factors through
  $\Sigma^{-n} \Sigma^{\infty} \Omega^{\infty} \Sigma^{n} X \ra X$,
  as one sees by suspending everything $n$ times and using that $\Sigma^\infty$
  is left adjoint to $\Omega^\infty$.
  Thus the composite $\Sigma^{-n} \Sigma^\infty W \ra X \ra Y$ is zero,
  and we have shown that $f : X \ra Y$ is a superphantom.
\end{proof}

Margolis states in his book~\cite[p.~81]{ma:ssa} that 
whether there exist non-zero superphantoms
is an open question.
We answer this question now.

\begin{prop}
  The mod $p$ \EM spectrum $\Hp$ is the source of a non-zero
  superphantom map.
\end{prop}

\begin{proof}
  We note that by the main result of~\cite{hora:ssh} there are
  no maps from $\Hp$ to a suspension spectrum.  Therefore,
  by the argument in the previous part of this section,
  $\Hp$ is not a retract of a wedge of suspension spectra.
  Therefore $\Hp$ is the source of a non-zero superphantom map.
\end{proof}

There are skeletal phantoms which are not superphantoms.
For example, following Gray~\cite{gr:ssn}
one can show that there are uncountably many skeletal phantoms
from $\CPi$ to $S^3$.  
But $\CPi$ is a suspension spectrum, and so none of these
maps is a superphantom.
We don't know if there is an example of a superphantom which is not a 
skeletal phantom.

\mysection{Topological ghosts} \label{se:ghosts}

In this section we continue to work in the stable homotopy category.
Again, there are three parts.
In the first, we describe the ghost projective class
and give its elementary properties.
In the second part, which is a little bit longer, but lots of fun, 
we calculate the length of $\RPn$ for small $n$.
And in the third part, we show that the Adams spectral sequence
with respect to the ghost projective class in a category
of \Ai modules over an \Ai ring is in fact a universal
coefficient spectral sequence, and we explain how this gives new
lower bounds on the ghost-length of a spectrum.

\subsection{The ghost projective class}

A map $X \ra Y$ is called a \dfn{ghost} if the induced map
$\pi_*(X) \ra \pi_*(Y)$ of homotopy groups is zero.
Let $\I$ denote the ideal of ghosts.
Let $\Proj$ denote the class of all retracts of wedges of spheres.
It is easy to see that $(\Proj,\I)$ is a projective class
(use Lemma~\ref{le:pc}) and that it generates (as $\Proj$ contains
the spheres).

Let's begin by noticing that there are spectra of arbitrarily
high length with respect to this projective class.
For example, the length of $\RP^{2^k}$ is at least $k+1$.
One sees this by noticing that if $u \in H^1(\RP^{2^k}; \Zt)$
is the non-zero class, then $\Sq^{2^{k-1}} \cdots \Sq^4 \Sq^2 \Sq^1 u$
is non-zero.  But the composite
\[ \RP^{2^k} \llra{u} \Sigma \Ht \xra{\Sq^1} \Sigma^2 \Ht 
   \xra{\Sq^2} \cdots \xra{\Sq^{2^{k-1}}} \Sigma^{2^k} \Ht \]
is in $\I^k$ and thus would be zero if $\RP^{2^k}$ had length $k$
or less.

On the other hand, by constructing $\RP^n$ one cell at a time,
it is clear that the length of $\RP^n$ is no more than $n$.
We will see in the next part of this section
that it is possible to improve on both of these
bounds.

The filtration of the morphisms of $\cS$ is also interesting.
Again by looking at composites of Steenrod operations, one sees
that the powers $\I^k$ are all non-trivial.
Also, every phantom map is a ghost.  The ghost-filtration of a
phantom map is analogous to what Gray called the ``index'' of a
phantom map, so we'll use that terminology here.
Every non-zero phantom map from the Moore spectrum $S(G)$ has index 1,
since $S(G)$ has length 2 with respect to the ghost ideal.
And we saw in Section~\ref{se:gms} that for $G$ non-zero and
divisible, such phantom maps exist.

The Lusternik--Schnirelmann category of a space $X$ is an upper bound
for the ``cup length'' of the reduced cohomology $H^*X$.  That is, if
the cup product $u_1 \cdots u_n$ is non-zero for some $u_i \in H^*X$,
then the Lusternik--Schnirelmann category of $X$ is at least $n$.
Stably there are no products in cohomology, but we have instead the
action of the Steenrod algebra.  And we saw above that if there is a
chain of Steenrod operations acting non-trivially on the mod 2 cohomology of
a spectrum $X$, say $\Sq^{i_1} \cdots \Sq^{i_n} u \neq 0$, then the
ghost-length of $X$ is at least $n$.  Thus we think of ghost-length as
a stable analogue of Lusternik--Schnirelmann category.

If for a pair $X$ and $Y$ of spectra the Adams spectral 
sequence abutting to $[X,Y]$ is strongly convergent, then
$\Io(X,Y)$ is zero.
Call a non-zero map in $\Io$ a \dfn{persistent ghost}.
We saw in Theorem~\ref{th:powers} that
$(\Projo,\Io)$ is a projective class,
but for all we know at this point, it could be that $\Projo$ contains
all of the objects of the stable homotopy category and that $\Io$
is zero.
Our first goal is to show that this is not the case.

\begin{prop} \label{pr:persistent}
  Let $X$ be a non-zero connective spectrum such that there are
  no maps from $X$ to a connective wedge of spheres.
  Then there is a persistent ghost $X \ra Y$ for some $Y$.
\end{prop}

If $X$ is a dissonant spectrum, such as $\Hp$, then there are no
maps from $X$ to a connective wedge of spheres.  
Indeed, a connective wedge of spheres is a (possibly desuspended) suspension
spectrum, and by a result of Hopkins and Ravenel~\cite{hora:ssh}
suspension spectra are harmonic.

\begin{proof}
  If $X$ is connective, the projectives $P_n$ in a ghost Adams 
  resolution for $X$ can be chosen to be connective wedges of spheres.  
  Let $W_n$ be the fibre of the map $X \ra X_n$.  
  There is a natural map $\Wedge W_{n} \ra X$ whose cofibre is
  a persistent ghost.
  So if there is no persistent ghost with $X$ as its source, then
  $X$ is a retract of $\Wedge W_n$.
  We saw at the beginning of Section~\ref{se:ass} that $W_{n}$
  lies in the cofibre sequence $W_{n-1} \ra W_n \ra P_{n-1}$.
  It follows inductively that if there are no maps from
  $X$ to a connective wedge of spheres, then
  there are no maps from $X$ to $W_n$ for each $n$.
  Thus the map $X \ra \Wedge W_n \ra \prod W_n$ is zero.
  But this map is also monic in homotopy groups, and so we conclude
  that $X$ is zero.
\end{proof}

Example~2.6 in the paper~\cite{li:atssshm} also implies that
$\Io $ is non-zero, but I have been unable to follow the argument.

Along the same lines, we can also obtain the next result.

\begin{prop} \label{pr:connective}
  If $X$ is a connective spectrum of length $n$, then $X$ can
  be built using $n$ connective wedges of spheres.
\end{prop}

We use the word ``built'' to mean ``built using cofibres and
retracts'', as in the definition of $\Proj_n$.

\begin{proof}
  Form an Adams resolution~\eqref{eq:ar} of $X$ with the $P_n$ chosen
  to be connective wedges of spheres.  
  Since $X$ has length $n$, the composite $X = X_0 \ra \cdots \ra X_n$ is
  zero.  Thus $X$ is a retract of $W_n$, the fibre of this composite.
  But we saw at the beginning of Section~\ref{se:ass} that
  $W_n$ can be built from $P_0, \ldots, P_{n-1}$,
  $n$ connective wedges of spheres.
\end{proof}

Similarly, one can show that a spectrum of finite type 
($\pi_i X$ a finitely generated abelian group for each $i$) 
and length $n$ can be built using $n$ wedges of spheres with only
a finite number of spheres of each dimension.

\begin{cor}
  If $X$ is a connective spectrum of finite length, then $X$ is harmonic.
\end{cor}

\begin{proof}
  We saw in the proof of Proposition~\ref{pr:persistent} that a 
  spectrum built from a finite number of
  connective wedges of spheres is harmonic, so the 
  result follows from Proposition~\ref{pr:connective}.
\end{proof}

\begin{question}
  Is every spectrum of finite length harmonic?  It would suffice to
  show that any wedge of spheres is harmonic.
  
  Another related question is whether every finite spectrum of 
  finite length can be built from a finite number of 
  \emph{finite} wedges of spheres.  We don't know the answer
  to either question.
\end{question}

\subsection{The ghost-length of real projective spaces} \label{se:rpn}

In this part of the section we give upper and lower bounds on the ghost-length
of $\RPn$.  The upper bound is obtained by building $\RPn$ carefully
using a cofibre sequence involving $\CPi$ and a Thom spectrum.
Our first lower bound is simply the length of the longest chain
of non-zero Steenrod operations acting on the mod 2 cohomology
of $\RPn$.  This bound agrees with the upper bound for $n < 20$, 
providing us with a calculation of the length of $\RP^n$ in this range.
However, we will show that while the squares 
$\Sq^1$, $\Sq^2$, $\Sq^4$ and $\Sq^8$ have ghost filtration 1, 
the squares $\Sq^{2^k}$ for $k \geq 4$ have ghost filtration
at least 2, and using this we obtain a significantly better lower
bound.  The ghost filtration of the Steenrod squares
is closely connected with the Hopf and Kervaire invariant problems,
and we give theorems explaining this relation.

We work localized at the prime 2.  The 2-local
category is triangulated, so all of our general theory applies.
We will write $H^{*} X$ for the mod 2 cohomology of $X$.

We begin by recalling the action of the Steenrod algebra on 
$H^*\RPi = \F_2[x]$, with $|x| = 1$.
The Steenrod square $\Sq^{2^k}$ acts non-trivially on $x^n$
if and only if the $k$th bit in the binary expansion of $n$ is 1.
Figure~\ref{fig:rpn} illustrates.

Now we describe a construction of $\RPi$ that was explained
to us by Mahowald.
The double cover map $S^1 \ra S^1$, with fibre $\Zt$, 
can be extended to a fibre sequence of spaces
\[ S^1 \llra{2} S^1 \lra \RPi \lra \CPi \llra{2} \CPi \]
by applying the classifying space functor.
Thus $\RPi$ is the circle bundle of the complex line bundle 
$\eta \tensor \eta$ over $\CPi$.  (Here $\eta$ denotes the
tautological line bundle.)  The Thom space is the
disk bundle modulo the circle bundle, and the disk bundle is homotopy
equivalent to $\CPi$, so there is a cofibre sequence
\[ \RPi \lra \CPi \lra \Th . \]
Writing $T$ for the desuspension of the Thom space,
we get a stable cofibre sequence
\[ T \lra \RPi \lra \CPi \]
which we will use to build $\RPn$ efficiently.
The Thom isomorphism tells us that $T$ can be built with a single 
cell in each odd non-negative dimension, and no other cells.
So $H^* T = H^* \Sigma \CPip$, where $\CPip $ denotes
$\CPi \Wedge S^{0}$.
In fact, the map $H^* T \la H^* \RPi$ is surjective,
so we can deduce the action of the Steenrod algebra on $H^* T$.
Figure~\ref{fig:rpn} displays the low degree part of the 
short exact sequence of modules over the Steenrod algebra that we obtain.
Each circle represents a basis element over $\F_2$ and the vertical
arrows give the action of the $\Sq^{2^k}$'s.
To make the pattern as nice as possible, we have replaced
$\RPi $ with $\RPip = \RPi \Wedge S^{0}$, and similarly with $\CPi $.
\newcommand{\sqi}{\ar[u]}
\newcommand{\sqiir}{\ar@/_.7pc/[uu]}
\newcommand{\sqiil}{\ar@/^.7pc/[uu]}
\newcommand{\sqiiir}{\ar@/_2.2pc/[uuuu]}
\newcommand{\sqiiil}{\ar@/^2.2pc/[uuuu]}
\newcount\d
\newcommand{\cell}{**[Fo:<5.5pt>] \d=12 \advance\d by-\Row \text{\footnotesize\the\d}}
{
\begin{figure}
\[ \label{eq:rpn} \xymatrixcolsep{3pc} \xymatrix{ \relax
  & \cell                & \cell \ar[l]                &           \\
  &                      & \cell                       & \cell \ar[l] \\
  & \cell                & \cell \ar[l] \sqi           &           \\
  &                      & \cell                       & \cell \ar[l] \\
  & \cell \sqiil \sqiiil & \cell \ar[l] \sqi \sqiil \sqiiil &         \\
  &                      & \cell \sqiir \sqiiir & \cell \ar[l] \sqiir \sqiiir\\
  & \cell \sqiiil        & \cell \ar[l] \sqi \sqiiil   &           \\
  &                      & \cell       \sqiiir        & \cell \ar[l] \sqiiir \\
  & \cell \sqiil         & \cell \ar[l] \sqi \sqiil    &           \\
  &                      & \cell        \sqiir         & \cell \ar[l] \sqiir \\
  & \cell                & \cell \ar[l] \sqi           &           \\
  &                      & \cell                       & \cell \ar[l]    \\
0 & H^* T \ar[l]         & H^* \RPip \ar[l]            & H^* \CPip \ar[l] 
                                                                     & 0 \ar[l]
} \]
\caption{The action of the Steenrod algebra}\label{fig:rpn}
\end{figure}
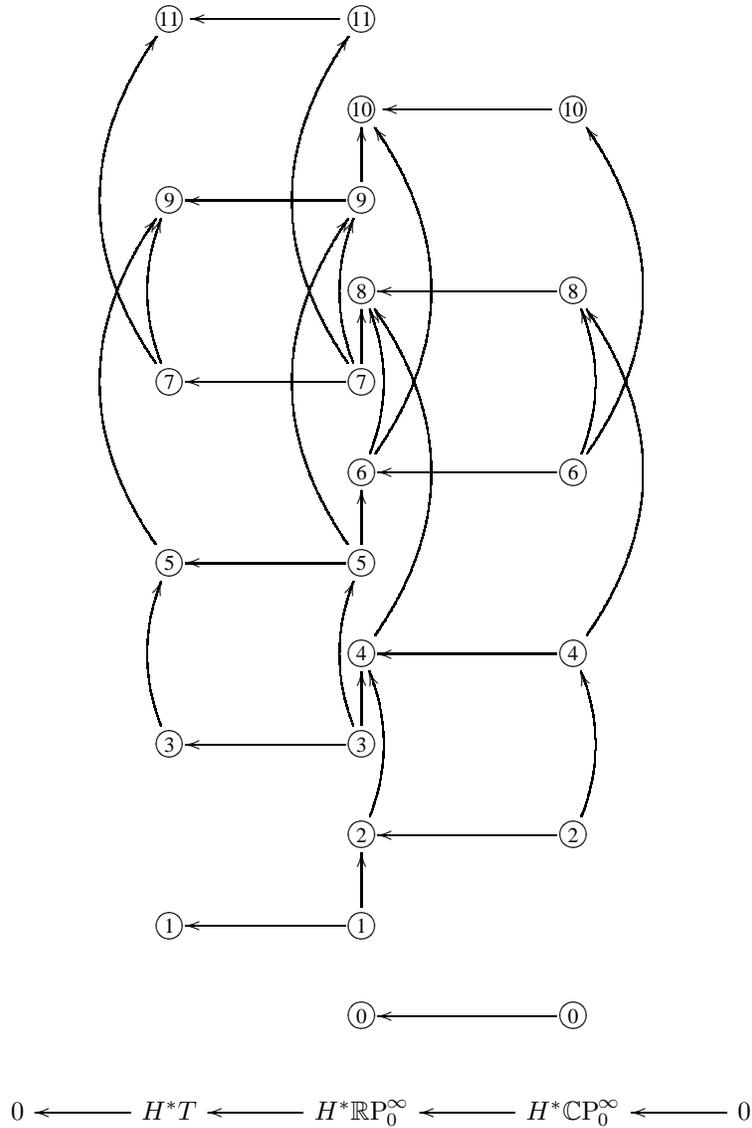
}

The essential point is that $\RPi$ can be built up odd cells first.
That is, we can first build $T$ completely, and then attach the even cells.
By looking at the cofibre sequence of skeleta, we see that
$\RPn$ can also be built by building up the odd part, and
then the even part.  We will use this to deduce an upper bound on
the length of $\RPn$. 

Let us begin by considering $\RP^3$.  The length of $\RP^3$ is at least
two, because $\Sq^1 x = x^2$.  And since $\RP^3$ has only three cells,
it has length at most three.  In fact, it has length two.  One way
to see this is to notice that the 3-skeleton of $T$
is $S^1 \Wedge S^3$, and so $\RP^3$ is formed by attaching a 2-cell
to this wedge.  

Both $\RP^4$ and $\RP^5$ can similarly be seen to have length three
because the 4- and 5-cells can to added to $\RP^3$ simultaneously.
(A connection between them would be detected by a non-zero $\Sq^1$
on $x^4$.)  But what about $\RP^6$?
This is trickier and will reveal the power of the decomposition
into odd and even cells.
We just mentioned that the 5-cell can be added after the 1-, 2- and
3-cells have been added.  But even cells are never needed for the
attachment of odd cells, so the 5-cell can actually be attached
at the \emph{same time} as the 2-cell.  And this means that the
6-cell can be attached at the same time as the 4-cell.  So
$\RP^6$ also has length three.  

With this under our belt, we now prove the following theorem.

\begin{thm}
  The length of $\RPn$ is no more than $\floor{n/4} + 2$. 
  Here $\floor{x}$ denotes the greatest integer less than or equal to $x$.
\end{thm}

\begin{proof}
  We prove this inductively by showing that we can add four cells at a time,
  if we are careful about the order.
  To make the pattern work from the start, we build $\RPnp = \RPn \Wedge S^{0}$
  instead of $\RPn$.  
  This makes little difference as both have the same length.
  We start with $S^1 \Wedge S^3$.
  To this we can add the 5-cell, since an odd cell only needs the odd
  cells below it.  Now in $T$ there is no $\Sq^2$ from
  the 5-cell to the 7-cell, and this implies that the 7-cell is not
  attached to the 5-cell.  So we can attach the 7-cell to $S^1 \Wedge S^3$.
  And we can of course attach the 0- and 2-cells.  
  Call the resulting complex $W$.
  $W$ has cells in dimensions 0, 1, 2, 3, 5 and 7, and has length 2.
  To $W$ we can attach the 9- and 11-cells, since they only require
  the odd cells below them and are not connected in $T$.
  At the same time we can attach the 4- and 6-cells, because they 
  are not connected in $\CPi$.  (Again, because there is no $\Sq^2$.)
  Thus we can add the 4-, 6-, 9- and 11-cells to $W$ to get a length
  3 complex $X$.
  In a similar way we see that we can add the 8-, 10-, 13- and 15-cells
  to $X$.  
  Thus $\RP^8$, $\RP^9$, $\RP^{10}$ and $\RP^{11}$ all have length at most 4.
  This pattern continues, proving the theorem.
\end{proof}

We saw in the previous part of this section that if there is a
chain of Steenrod operations acting non-trivially on the cohomology of
a spectrum $X$, say $\Sq^{i_1} \cdots \Sq^{i_n} u \neq 0$, then the
ghost-length of $X$ is at least $n$.  
Letting $\Stl(\RPnp)$ denote one more than 
the length of the longest such chain in the 
cohomology of $\RPnp$, one obtains the following sequence of numbers,
\newcommand{\nts}{\negthickspace}
{\small
\[ \label{eq:stl}\nts\begin{array}{c|*{22}{>{\PBS{\raggedleft}\hspace{0pt}}p{1.5em}@{}}}
n & \nts -1 & 0 & 1 & 2 & 3 & 4 & 5 & 6 & 7 & 8 & 9 & 10 & 11 & 12 & 13
                            & 14 & 15 & 16 & 17 & 18 & 19 & 20 \\ 
  \hline
\nts\Stl(\RPnp) & \nts 0 & 1 & 1 & 2 & 2 & 3 & 3 & 3 & 3 & 4 & 4 &  4 &  4 &  5 &  5 
                            &  5 &  5 &  6 &  6 &  6 & 6 & 6 \comma 
\end{array} \]
}%
where we regard $\RP^{-1}_{0}$ as the zero object, which has
length zero.
If we count the number of consecutive 0's, then the number of
consecutive 1's, and so on, we obtain the sequence 1, 2, 2, 4,
4, 4, 8, 8, 8, 8, $\ldots$.  

\begin{thm}[Vakil~\cite{va:slrps}] \label{th:ravi}
  The sequence obtained in this way consists of the powers of\/ $2$, in
  order, with $k+1$ repetitions of\/ $2^k$.              \qed
\end{thm}

This theorem completely determines the sequence $\Stl(\RPnp)$.
The proof of the theorem has a  striking feature. 
Vakil studies a more fundamental sequence defined in the following way.  
The $n$th term of the sequence is one more than
the length of the longest chain of non-zero Steenrod operations
in $\RPi$ which \emph{ends at the $n$th cell}.
Starting with $n = 1$, this sequence begins 1, 2, 1, 3, 2, 3, 1, $\ldots$.
The $n$th term of the sequence displayed in the table 
\vpageref[above]{eq:stl} is obtained by taking the
supremum of the first $n$ terms of the fundamental sequence.
The method of proof that Vakil uses to determine where the
jumps in the values of the suprema occur is to explicitly
define certain ``canonical moves''.
In more detail, given an $n$, Vakil determines a non-zero Steenrod
operation $\Sq^{2^{k}}$ which is the
last step in a longest chain ending at the $n$-cell 
(if there are any Steenrod operations hitting the $n$-cell).  
Using this, one can quickly compute the $n$th term in the fundamental
sequence by following the canonical moves downwards until one reaches
a cell not hit by a Steenrod operation.
And one can also use the canonical moves to prove Theorem~\ref{th:ravi}.

\begin{note}
  For $2 \leq n \leq 19$, $\Stl(\RPn) = \floor{n/4} + 2$.  
  Thus we know the length of $\RPn$ for such $n$.
  But $\Stl (\RP^{20} ) = 6$ and $\floor{20/4} + 2 = 7$,
  and for larger $n$ this just get worse.
  For example, $\Stl (\RP^{2^{20}}) = 136$
  and $\floor{2^{20}/4} + 2 = 2^{18} + 2$.
\end{note}

One might wonder whether the lower bound is correct.  It is not.  For
example, the length of $\RP^{2^{20}}$ is actually at least 264.
The first case where I know that the Steenrod length gives the
wrong answer is for $n = 56$.  The Steenrod length of $\RP^{56}$
is 10 but the length is at least 11.
These facts are deduced from the following result.

\begin{thm}\label{th:hio}
  The Steenrod operations $\Sq^1$, $\Sq^2$, $\Sq^4$ and $\Sq^8$ have
  ghost-filtration exactly one.
  The Steenrod operations $\Sq^{16}$, $\Sq^{32}$, $\ldots$ have
  ghost-filtration at least two.
\end{thm}

This allows us to count a higher Steenrod square occurring in the cohomology
of $\RPn$ as two maps.  Using a computer to do the computation,
this is how we obtained the improved lower bound on the length
of $\RP^{2^{20}}$.  
It turns out that the improved lower bound is still wrong in general.
For example, the improved lower bound tells us that the
length of $\RP^{127}$ is at least 17.
But it also tells us that the length of
$\RP^{128}$ is at least 19, hence
the length of $\RP^{127}$ must be at least 18.

\begin{proof}[Proof of Theorem~\ref{th:hio}]
  The (non-identity) Steenrod squares all have filtration at least one.
  It is well-known that $\Sq^1$, $\Sq^2$, $\Sq^4$ and $\Sq^8$ act 
  non-trivially on the length two complexes
  $\RP^2$, $\CP^2$, $\HP^2$ and $\OP^2$ respectively,
  so these operations must have filtration exactly one.

  To show that the higher squares can be factored into two pieces,
  each zero in homotopy, we make use of Adams' result on the
  Hopf invariant one problem~\cite{ad:neehio}.
  Adams shows that no complex with only two cells supports a
  non-zero $\Sq^{2^k}$ with $k \geq 4$.
  Consider the beginning of a ghost Adams resolution of $H$, the
  mod 2 \EM spectrum:
\[\xymatrixcolsep{0.6pc} 
  \xymatrixrowsep{1.05pc}
  \xymatrix{ \relax
   H \ar[rr] && {\Bar{H}} \ar[ld] \circar[ld] \ar[rr] && {\Bar{\Bar{H}}} \ar[ld] \circar[ld] && \\ %
   & S^{0} \ar[lu] && P \ar[lu] && \period %
} \]
  Here $P$ is a large wedge of spheres.
  The fibre $W$ of the map $H \ra \Bar{\Bar{H}}$ has length two---it lies
  in a cofibre sequence $S^{0} \ra W \ra P$.
  If the composite 
\[
W \lra H \xra{\Sq^{2^{k}}} \Sigma^{2^{k}} H
\]
  is zero, then $\Sq^{2^{k}}$ factors through 
  $H \ra \Bar{H} \ra \Bar{\Bar{H}}$ and thus has filtration at least two.
  So we have reduced the problem from checking that $\Sq^{2^k}$
  vanishes on \emph{all} cohomology classes of \emph{all} length two
  spectra to checking that it is zero on a \emph{particular} cohomology 
  class in a \emph{particular} length two spectrum.
  To do this, notice that the composite $W \ra H \ra \Sigma^{2^{k}} H$
  factors (uniquely) through the map $W \ra P$.  To show that the map
  $P \ra \Sigma^{2^{k}} H$ is zero, it suffices to check this
  on each $2^{k}$-dimensional sphere appearing as a summand of $P$.
  Choose such a summand, and consider the following diagram
  \[ \xymatrix{ \relax
  & \Sigma^{2^{k}} H \\
  &                H \ar[u]^{\Sq^{2^k}} \\
  S^{0} \ar[r] & W \ar[u] \ar[r] & P \ar[uul] \ar[r] & S^{1} \\
  S^{0} \ar@{=}[u] \ar[r] & W' \ar[u] \ar[r] & S^{2^k} \ar[u] \ar[r] & S^{1} \ar@{=}[u] \comma
  } \]
  in which $W'$ is defined to be the fibre of the map from $S^{2^{k}}$
  to $S^{1}$ and the map $W' \ra W$ is some choice of fill-in map.
  We must show that the composite from $S^{2^{k}}$ to $\Sigma^{2^k} H$
  is zero.  Well, by Adams' result, it is zero when restricted to
  $W'$; so it factors through $S^{1}$;  but there are no maps
  from $S^{1}$ to $\Sigma^{2^{k}} H$; so it must be null.
\end{proof}

Now we quote a theorem which relates the filtration of the
Steenrod squares to the Kervaire invariant problem.

\begin{thm}[W.-H. Lin~\cite{li:dass}]\label{th:lin}
If, for $k \geq 4$, the Kervaire class $\theta_{k-1}$ exists
and has order 2, then $\Sq^{2^{k}}$ has filtration exactly 2.  \qed
\end{thm}

We end this section with a conjecture.

\begin{cnj}
The length of $\RPn$ is an increasing function of $n$.\footnote{This has
now been proven by the author.}
\end{cnj}

\subsection{A universal coefficient spectral sequence}\label{sse:ucss}

In this part of the section we need to briefly step outside of
the homotopy category.
Given a ring spectrum $R$, we would like to have a triangulated
category of $R$-modules.
Unfortunately, this isn't possible if $R$ is simply
a monoid object in the homotopy category.
So by an ``\Ai ring spectrum'' $R$ we mean any notion of
structured ring spectrum such that the homotopy category $\RMod$ 
of ``\Ai module spectra'' is triangulated and satisfies the following formal
properties.
There is a ``free module'' functor $F : \cS \ra \RMod$
which is left adjoint to a ``forgetful'' functor $U : \RMod \ra \cS$.
Both $F$ and $U$ preserve triangles, commute with suspension,
and commute with coproducts,
and the composite $U F$ is naturally isomorphic to the functor
sending $X$ to $R \Smash X$.
We will usually omit writing $U$, and will write $R \Smash X$
for both $F X$ and $U F X$, with the context making clear which
is intended.

There are various notions of structured ring spectra available to us 
today~\cite{ekmm:rmasht,hoshsm:ss,lemast:esht}.  Unfortunately, 
we know of no published proof that the formal properties hold
in these settings.  It is certainly expected that they do.

Fix an \Ai ring spectrum $R$ and
write $R$ for $F S^{0}$.  $R$ is the ``sphere'' in the category
of \Ai $R$-modules.
Indeed, by adjointness, $[R,M]_{R} = [S^{0},M] = \pi_{0}M$,
where we write $[M,N]_{R}$ for maps from $M$ to $N$ in \RMod.

Because $F$ preserves triangles, commutes with suspension, and
preserves retracts, it is clear that if a spectrum $X$ can be built
from $n$ wedges of spheres, then $F X$ can be built from $n$ wedges
of suspensions of $R$.
To make this more precise, we note that in \RMod there is a 
projective class $(\Proj_{R},\I_{R})$, where  $\Proj_{R}$ 
is the collection of retracts of wedges of suspensions of $R$
and $\I_{R}$ is the collection of maps zero in homotopy groups.
To put it another way, $\Proj_{R}$ is the image of $\Proj$
under $F$ (with retracts thrown in), and $\I_{R}$ is $U^{-1} \I$.
And our claim is that the length of $F X$ with respect to
$(\Proj_{R},\I_{R})$ is no more than the length of $X$ with
respect to $(\Proj,\I)$.

So it would be useful to give a lower bound for the length
of an \Ai $R$-module.  This is accomplished in the remainder of
the section.

\begin{thm}\label{th:RModSS}
Let $M$ and $N$ be \Ai modules over the \Ai ring spectrum $R$.
Then there is a conditionally convergent spectral sequence
\[
E_{2}^{*,*} = \Ext_{R_{*}}^{*,*}(M_{*},N_{*}) \implies [M,N]_{R}^{*} .
\]
If $M$ has length at most $n$ with respect to $(\Proj_{R},\I_{R})$,
then $E_{n+1} = E_{\infty}$.
\end{thm}

By taking $M = R \Smash X$ we get the following consequence.

\begin{cor}\label{co:UCSS}
If $X$ is a spectrum and $N$ is an \Ai module over an \Ai ring spectrum $R$,
then there is a conditionally convergent spectral sequence
\[
E_{2}^{*,*} = \Ext_{R_{*}}^{*,*}(R_{*}X,N_{*}) \implies N^{*}X .
\]
If $X$ has length at most $n$ with respect to $(\Proj,\I)$,
then $E_{n+1} = E_{\infty}$.
\qed
\end{cor}

This is called the \dfn{universal coefficient spectral sequence}.
For another account, see~\cite{elgrkrma:cashtct}.

\begin{proof}[Proof of Theorem]
The spectral sequence is simply the Adams spectral sequence 
with respect to the projective class $(\Proj_{R},\I_{R})$.
The $E_{2}$-term consists of the derived functors of $[-,N]_{R}$
applied to $M$.
A projective resolution of $M$ is a sequence
\[
0 \lla M \lla P_{0} \lla P_{1} \lla \cdots
\]
of \Ai $R$-modules with each $P_{s}$ in $\Proj_{R}$ and
which is exact in homotopy.
Thus applying $\pi_{*}(-)$ gives a projective resolution of $M_{*}$.
One can check that $[P_{s},N]_{R} = \Hom_{R_{*}}(\pi_{*}P_{s},N_{*})$
and thus that the $E_{2}$-term is $\Ext_{R_{*}}^{*,*}(M_{*},N_{*})$.
The projective class generates and $\I_R$ is closed under coproducts,
so by Proposition~\ref{pr:condconv}, the spectral sequence is 
conditionally convergent.

That the spectral sequence collapses at $E_{n+1}$ when $M$ has length
at most $n$ is Proposition~\ref{pr:collapse}.
\end{proof}

Thus the existence of a non-zero differential $d_{n}$ implies that
$X$ has length at least $n$.
We suspect that for $X = \RPn$, $R = J$ or $KO$, and $N = KO$,
this gives a very good lower bound for the ghost-length of $\RPn$.
However, while in some cases we have been able to compute the 
$E_{2}$- and $E_{\infty}$-terms,
we haven't been able to conclude anything about the differentials.

We finish this section by mentioning the following example
of Theorem~\ref{th:RModSS}.
Take $R = S^{0}$ and $M = N = H$, the mod 2 \EM spectrum.
Then we get a spectral sequence with $E_{2}$-term
\[
E_{2}^{*,*} = \Ext_{\pi_{*}S^{0}}^{*,*}(\Zt,\Zt)
\]
converging to the Steenrod algebra.
This has been called the dual or reverse Adams spectral sequence and
has been studied by various authors~\cite{li:dass,li:atssshm}.
An unstable version is described in~\cite{bl:hssh} and~\cite{bl:orrass}.

\mysection{Algebraic ghosts} \label{se:alg-ghosts}

We now discuss a projective class which provided the motivation for
this work.
Indeed, an old result of Kelly~\cite{ke:cmizhm} (presented here as
Theorem~\ref{th:kelly}) concludes that under certain conditions a 
composite of maps vanishes.  We wondered whether there was more
than just a superficial similarity between this result and the
fact that in the stable homotopy category a composite of two
phantoms is zero.
It turns out that the arguments can be arranged to have a common part,
by proving that the ideals in question are parts of projective
classes and then applying Theorem~\ref{th:powers}.

We work in the derived category of an abelian category in this section,
and so we begin with a brief overview of the derived category.
Good references are~\cite{we:iha} and~\cite{iv:cs}.

Let $\cA$ be an abelian category with enough projectives.
We mean this in the usual sense---that is, we are assuming
that the categorical projectives and the categorical epimorphisms
form a projective class.
We also assume that $\cA$ satisfies Grothendieck's AB 5 axiom which says that
set-indexed colimits exist and filtered colimits are exact~\cite{gr:qpa}.  
We write $\Ch$ for the category of $\Z$-graded chain complexes
of objects of $\cA$ and degree 0 chain maps.
To fix notation, assume the differentials have degree $-1$.
For an object $X$ of $\Ch$, define $Z_n X := \ker(d : X_n \ra X_{n-1})$ and
$B_n X := \im(d : X_{n+1} \ra X_n)$, and write $H_n X$ for the quotient.

Write $\K$ for the category in which we identify chain homotopic maps.
It is well-known that the category $\K$ is triangulated, so
we will only briefly describe the triangulation.
There is an automorphism $\Sigma$ of $\Ch$ which is defined
on objects by $(\Sigma X)_n = X_{n-1}$ and $d_{\Sigma X} = -d_X$,
and on morphisms by $(\Sigma f)_n = f_{n-1}$;
this induces an automorphism $\Sigma$ of $\K$ which serves
as the suspension for the triangulated structure.
A short exact sequence
\[ 0 \lra W \llra{i} X \llra{p} Y \lra 0 \]
of chain complexes is \dfn{weakly split} if for each $n$ the
sequence
\[ 0 \lra W_n \llra{i_n} X_n \llra{p_n} Y_n \lra 0 \]
is split.
Given a weakly split short exact sequence of chain complexes as above,
choose for each $n$ a splitting of the $n$th level, \ie choose
maps $q_n : X_n \ra W_n$ and $j_n : Y_n \ra X_n$ such that
$p_n j_n = 1$, $q_n i_n = 1$ and $i_n q_n + j_n p_n = 1$.
Define $h_n : Y_n \ra W_{n-1}$ to be $q_{n-1} d_X j_n$.
One easily checks that $h$ is a chain map $Y \ra \Sigma W$
and that up to homotopy $h$ is independent of the choice of
splittings.
A triangle in $\K$ is a sequence isomorphic (in $\K$) to one of the
form $W \ra X \ra Y \ra \Sigma W$ constructed in this way 
from a weakly split short exact sequence.
See~\cite{iv:cs} for details.

One fact we use about the triangulation is that
the homology functors $H_n : \K \ra \cA$ send triangles to
long exact sequences.  

A chain map $f : X \ra Y$ is a \dfn{quasi-isomorphism} if it
induces an isomorphism in homology.
The derived category $\D$ is the category obtained from $\Ch$ 
by formally inverting the quasi-isomorphisms.
(It is equivalent to invert the quasi-isomorphisms in $\K$.)
With our hypotheses on $\cA$, this category of fractions
exists~\cite[Exercise~10.4.5]{we:iha}.
In fact, it is equivalent to the full subcategory of $\K$
containing the ``cofibrant'' complexes.
A complex $X$ is \dfn{cofibrant} if it can be written as
an increasing union $X = \cup_{n \geq 0} C^n$ of subcomplexes $C^n$
with $C^0 = 0$ and $C^n/C^{n-1}$ a complex of projectives
with zero differential.
To prove the equivalence of categories, 
one shows that for any $X$ there is a cofibrant complex $W$
and a quasi-isomorphism $W \ra X$,
and that when $X$ is cofibrant, the natural
map $\K(X,Y) \ra \D(X,Y)$ is an isomorphism for all $Y$.

The derived category is a triangulated category.
The automorphism $\Sigma$ of $\Ch$ 
induces an automorphism $\Sigma$ of $\D$.
There is a natural functor $\K \ra \D$, and
a sequence $X \ra Y \ra Z \ra \Sigma X$ is a triangle in $\D$ if and only
if it is isomorphic (in $\D$) to the image of a triangle in $\K$.
One important fact about the triangulation is that 
if $f : X \ra Y$ is a chain map which is an epimorphism in
each degree, then the fibre of $f$ in $\D$ is given by
the degreewise kernel.

We record the following lemma, whose proof is straightforward.

\begin{lemma} \label{le:proj-equiv}
  Let $X$ be an object of $\K$.  Then the following are equivalent:
  \begin{roenumerate}
  \item $X$ is isomorphic in $\K$ to a complex of projectives with
    zero differential.
  \item $X$ is isomorphic in $\K$ to a complex $Y$ with $Y_n$, $Z_n Y$, 
    $B_n Y$ and $H_n Y$ projective for each $n$.
  \item $X$ is isomorphic in $\K$ to a complex $Y$ with 
    $B_n Y$ and $H_n Y$ projective for each $n$.                   \qed
  \end{roenumerate}
\end{lemma}

Since the homology of a complex is analogous to the homotopy of
a spectrum, we call a map which is zero in homology a \dfn{ghost}.
Let $\I$ denote the ideal of ghosts in $\D$.
Call an object $P$ of $\D$ \dfn{ghost projective} if it is
isomorphic (in $\D$) to an object satisfying the equivalent conditions
of Lemma~\ref{le:proj-equiv}.
Write $\Proj$ for the collection of ghost projective complexes.
One can check that $\Proj$ is closed under retracts.

As the reader has no doubt guessed, we have the following result.

\begin{prop} \label{pr:algh}
  The pair $(\Proj,\I)$ forms a projective class.
\end{prop}

\begin{proof}
  We begin by showing that $\Proj$ and $\I$ are orthogonal.
  Let $P$ be a ghost projective complex.  
  Without loss of generality, we may assume that $P$ is a 
  complex of projectives with zero differential.
  Since a complex with zero differential
  is a coproduct of complexes concentrated in a single degree,
  we may even assume that $P$ is a 
  projective object concentrated in degree zero, say.
  Such a complex is cofibrant, so $\D(P,Y) = \K(P,Y)$ for any $Y$.
  Now suppose that $f : P \ra Y$ is a map.
  That is, we have a map $f : P_0 \ra Y_0$ such that the composite
  $P_0 \ra Y_0 \ra Y_{-1}$ is zero; so $f$ factors through
  the kernel to give a map $P \ra Z_0 Y$.
  If $f$ is zero in homology, then the composite
  $P_0 \ra Z_0 Y \ra H_0 Y$ is zero; so $f$ factors through
  the inclusion of $B_0 Y$ into $Z_0 Y$.
  And because $P_0$ is projective, $f$ lifts over the 
  epimorphism $Y_1 \ra B_0 Y$.  
  That is, $f$ is null homotopic.
  We conclude that if $P$ is ghost projective and $g : P \ra X$ and
  $h : X \ra Y$ are maps in $\D$ with $h$ zero in homology,
  then the composite is zero in $\D$.

  Now, given a chain complex $X$, we construct a cofibre sequence
  $P \ra X \ra Y$ with $P$ ghost projective and with $X \ra Y$ zero in
  homology.
  First we choose projectives $P^{B_n}$ and $P^{H_n}$ and epimorphisms
  $P^{B_n} \ra B_n X$ and $P^{H_n} \ra H_n X$.   
  It is easy to see that we can now choose a projective $P^{Z_n}$
  and an epimorphism $P^{Z_n} \ra Z_n X$ which fit into a
  diagram
  \[ \xymatrix{
     0 \ar[r] & 
     P^{B_n} \ar[r] \ar@{{}->>}[d] & 
     P^{Z_n} \ar[r] \ar@{{}->>}[d] & 
     P^{H_n} \ar[r] \ar@{{}->>}[d] & 
     0 \\
     0 \ar[r] &
     B_n X \ar[r] &
     Z_n X \ar[r] &
     H_n X \ar[r] &
     0
  } \]
  with exact rows.
  Similarly, one can form a diagram
  \[ \xymatrix{
     0 \ar[r] & 
     P^{Z_n} \ar[r] \ar@{{}->>}[d] & 
     P^{X_n} \ar[r] \ar@{{}->>}[d] & 
     P^{B_n} \ar[r] \ar@{{}->>}[d] & 
     0 \\
     0 \ar[r] &
     Z_n X \ar[r] &
     X_n \ar[r] &
     B_n X \ar[r] &
     0
  } \]
  with exact rows and with $P^{X_n}$ projective.
  Defining $P_n := P^{X_n}$ and using the composite
  $P^{X_{n+1}} \ra P^{B_n} \ra P^{Z_n} \ra P^{X_n}$ as a differential,
  we get a chain complex $P$.
  By definition, $Z_n P = P^{Z_n}$.  The same holds for $B_n$ and $H_n$,
  so, by Lemma~\ref{le:proj-equiv}, $P$ is ghost projective.
  The maps $P^{X_n} \ra X_n$ piece together to give a chain map $P \ra X$.
  Under the functor $H_n$ this map induces the chosen epimorphism
  $P^{H_n} \ra H_n X$, and since $H_n$ is an exact functor, 
  the cofibre map $X \ra Y$ is zero in homology.

  Thus by Lemma~\ref{le:pc} we have a projective class.
\end{proof}

The main result of this section is the following theorem.

\begin{thm} \label{th:algh}
  Let $X$ be a complex such that the projective dimensions of
  $B_n X$ and $H_n X$ are less than $k$ for each $n$.
  Then the projective dimension of $X$ 
  with respect to the ideal of ghosts
  is less than $k$.
  In particular, $X$ has length at most $k$, and a $k$-fold composite
  \[ X \lra Y^1 \lra \cdots \lra Y^k \]
  of maps each zero in homology is zero in $\D$.
\end{thm}

\begin{proof}
  Let $X^0 = X$.
  As in the proof of the previous proposition, one can construct a
  map $P^0 \ra X^0$ such that each of the
  maps $P^0_n \ra X^0_n$, $H_n P^0 \ra H_n X^0$, $B_n P^0 \ra B_n X^0$
  and $Z_n P^0 \ra Z_n X^0$ is an epimorphism from a projective.
  Let $X^1$ be the suspension of the degreewise kernel, which is a choice 
  of cofibre, and inductively continue this process.
  For any exact sequence
  \[ 0 \lra A \lra Q_{k-2} \lra \cdots \lra Q_1 \lra H_n X \lra 0 \]
  in $\cA$ with each $Q_i$ projective, 
  the object $A$ is projective because of the assumption on $H_n X$.  
  The same holds with $H_n X$ replaced with $B_n X$, $Z_n X$ and $X_n$,
  for each $n$, so applying Lemma~\ref{le:proj-equiv} 
  one finds that $X^{k-1}$ is ghost projective.
  Thus $X^{k-1}$ has length at most one, $X^{k-2}$ has length at 
  most two, and inductively, $X = X^0$ has length at most $k$.
\end{proof}

\begin{cor}\label{co:gpd}
  If every object in $\cA$ has projective dimension less than $k$,
  then every object of $\D$ has projective dimension less than $k$ with respect
  to the ideal of ghosts.                                            \qed
\end{cor}

Note that the projective dimension of $H_{n} X$ is a lower bound
for the projective dimension of $X$.  

By assuming that $X$ is a complex of projectives, we can strengthen 
the conclusion of the theorem and obtain
the following result of Kelly~\cite{ke:cmizhm}.

\begin{thm} \label{th:kelly}
  Let $X$ be a complex of projectives such that 	
  the projective dimensions of $B_n X$ and $H_n X$ are less than $k$
  for each $n$.
  Then the projective dimension of $X$ 
  with respect to the ideal of ghosts in $\K$
  is less than $k$.
  In particular, $X$ has length at most $k$, and a composite
  \[ X \lra Y^1 \lra \cdots \lra Y^k \]
  of $k$ maps in $\Ch$, each zero in homology, is null homotopic.
\end{thm}

We emphasize that we are claiming that the composite is null
homotopic, not just zero in the derived category.

\begin{proof}[Sketch of proof]
One begins by showing that the collection
of retracts (in $\K$) of complexes satisfying the conditions of
Lemma~\ref{le:proj-equiv} along with the ideal of maps in $\K$
which are zero in homology is a projective class.  Then one
imitates the proof of Theorem~\ref{th:algh}, making use of the fact that
if $Z$ is a complex of projectives and
$Y \ra Z$ is a map in $\Ch$ which is degreewise surjective,
then the complex $X$ of degreewise kernels is the fibre
(since the sequence $X \ra Y \ra Z$ is degreewise split).
\end{proof}

\begin{cor}\label{co:hocof}
  Let $X$ be a complex of projectives such that 	
  the projective dimensions of $B_n X$ and $H_n X$ are less than $k$
  for each $n$.
  Then $X$ has the homotopy type of a cofibrant complex.
  That is, $X$ is isomorphic in $\K$ to a cofibrant complex.
\end{cor}

\begin{proof}
  This follows from Theorem~\ref{th:kelly} and the following lemma.
\end{proof}

\begin{lemma}
  A complex in $\K$ of finite length has the homotopy type of
  a cofibrant complex.
\end{lemma}

\begin{proof}
A complex $X$ has the homotopy type of a cofibrant complex
if and only if the natural map $\K(X,Y) \ra \D(X,Y)$ is an
isomorphism for all $Y$.
A complex of projectives with zero differential is cofibrant,
and so a retract in $\K$ of such a complex
has the homotopy type of a cofibrant complex.
The functors $\K(-,Y)$ and $\D(-,Y)$ are exact and send coproducts
to products, so the collection of complexes of the homotopy type
of a cofibrant complex is closed under coproducts and cofibre sequences.
Thus this collection contains all complexes of finite length.
\end{proof}

One can also prove Corollary~\ref{co:hocof} directly and
then deduce Theorem~\ref{th:kelly} from Theorem~\ref{th:algh}.

\mysection{Algebraic phantom maps} \label{se:alg-phantoms}

In this section we study phantom maps in the derived category of
an associative ring $R$.
We restrict attention from a general abelian category to the
category of $R$-modules because it is in this setting that one
can easily discuss the notion of purity.
We provide such a discussion in the first part of this section.
In the second part we describe the finite objects in the derived
category of $R$ and the phantom projective class that results.
Under some assumptions on $R$ we show that there is
a relation between pure extensions and phantom maps.
We end by recounting an example of Neeman's that shows
that phantom maps can compose non-trivially and hence that
Brown representability can fail in the derived category of
$R$-modules.

In this section $\Ch$ denotes the category of chain complexes of $R$-modules, 
$\K$ denotes the category obtained from $\Ch$ 
by identifying chain homotopic maps, 
and $\D$ denotes the derived category obtained from $\Ch$
by inverting quasi-isomorphisms.
See Section~\ref{se:alg-ghosts} for descriptions of these
categories.

We write $\Hom$ for $\Hom_R$ and $\tensor$ for $\tensor_{R}$, and,
unless otherwise stated, we take our modules to be left $R$-modules.

\subsection{Purity}

A module $P$ is \dfn{pure projective} if it is a summand of a coproduct
of finitely presented modules.
A sequence $K \ra L \ra M$ is \dfn{pure exact} if it is exact 
under $\Hom(P,-)$ for each pure projective $P$.
A longer sequence is \dfn{pure exact} if each three term subsequence
is pure exact.
A map $L \ra M$ is a \dfn{pure epimorphism} if each map $P \ra M$
from a pure projective factors through $L$.

We recall some standard facts about purity.
A good reference here is~\cite{wa:pacm}.

\begin{prop}\label{pr:enough-pps}
\begin{roenumerate}
\item The pure projectives, pure exact sequences and pure epimorphisms
      form a projective class as described in Sections~\ref{sse:pointed} 
      and~\ref{sse:wk}.
\item Every projective is pure projective;
      every pure exact sequence is exact;
      and every pure epimorphism is epic.
\item An infinite sequence
      \[ \cdots \lra M_{-1} \lra M_{0} \lra M_{1} \lra \cdots \]
      is pure exact if and only if it is exact after tensoring with
      each right module $E$.
      In particular, a finite sequence
      \[ 0 \lra M_{0} \lra M_{1} \lra \cdots \lra M_{k} \lra 0 \]
      beginning and ending with $0$ is pure exact if and only if it
      is exact after tensoring with each right module $E$.
      A similar statement holds for semi-infinite sequences
      beginning or ending with $0$.                          \qed
\end{roenumerate}
\end{prop}

\begin{note}
The sequence $0 \lra \Z \llra{2} \Z$ of abelian groups is pure exact,
but fails to be exact after tensoring with $\Zt$.
\end{note}

For $k \geq 1$, a \dfn{pure extension of length $k$} is a pure exact sequence
\[ 0 \lra N \lra A_1 \lra \cdots \lra A_k \lra M \lra 0 . \]
A morphism of extensions is a commutative diagram of the form
  \[ \xymatrix{
  E: & 0 \ar[r] & N \ar[r] \ar@{=}[d] & A_1  \ar[r] \ar[d] & 
          \cdots \ar[r] & A_k  \ar[r] \ar[d] & M \ar[r] \ar@{=}[d] & 0 \\
  E':& 0 \ar[r] & N \ar[r] & A_1' \ar[r]  & 
          \cdots \ar[r] & A_k' \ar[r] & M \ar[r] & 0 \period
  } \]
We say that two pure extensions $E$ and $E'$ are equivalent 
if they are connected by a chain of morphisms of pure extensions,
with the morphisms going in either direction.
Under the operation of Baer sum, the collection
$\PExt^k(M,N)$ of equivalence classes of pure
extensions of length $k$ forms an abelian group
which is a functor of $M$ and $N$:
the induced maps are given by pullback and pushforward of extensions.
The functors $\Ext^{k}$ can be defined in the same way.
(See~\cite{ma:h} for details.)
By forgetting pure exactness, one gets
a natural transformation $\PExt^k \ra \Ext^k$
which is not always a monomorphism.
We set $\PExt^0(M,N) = \Ext^{0}(M,N) = \Hom(M,N)$.

Extensions $\alpha \in \PExt^{k}(K,L)$ and $\beta \in \PExt^{l}(L,M)$
can be spliced together to give their \dfn{Yoneda product}, an
element of $\PExt^{k+l}(K,M)$ which we denote $\beta \alpha$.
We also use this composition notation when one or both of
$k$ and $l$ are zero.
Similarly, one can compose extensions in $\Ext$, and the natural
transformation $\PExt \ra \Ext$ respects composition.

As one might expect, $\PExt^*(M,N)$ can be calculated by forming a 
pure projective resolution of $M$, applying $\Hom(-,N)$, 
and taking homology.
In fact, associated to each filtered diagram $\{M_{\alpha}\}$ of 
finitely presented modules with colimit $M$ there is a natural 
pure projective resolution of $M$.  
(That every $M$ is in fact a filtered colimit of 
finitely presented modules is proved below.)
Consider the following sequence:
\begin{equation}\label{eq:ppr} 
   \cdots \lra
   \bigoplus_{\alpha \ra \beta \ra \gamma} M_\alpha \lra
   \bigoplus_{\alpha \ra \beta} M_\alpha \lra
   \bigoplus_{\alpha} M_\alpha \lra 
   M \lra 
   0 . 
\end{equation}
The sums are over sequences of morphisms in the filtered diagram.
Write $i_{\alpha} : M_{\alpha} \ra M$ for the colimiting cone to $M$.
The map $\oplus_\alpha M_\alpha \ra M$ is equal to $i_{\alpha}$
on the $\alpha$ summand.
A summand of $\oplus_{\alpha \ra \beta} M_\alpha$ is indexed by a
triple $(\alpha,\beta,u)$, where $u$ is
a map $M_\alpha \ra M_\beta$ such that $i_\alpha = i_\beta u$.
The map $\oplus_{\alpha \ra \beta} M_\alpha \ra \oplus_\alpha M_\alpha$
sends the summand $M_\alpha$ indexed by such a triple to the
$M_\alpha$ summand using the identity map and to the $M_\beta$ summand
using the map $-u$.
In general, one gets an alternating sum.
Taking cohomology gives the derived functors of colimit
(see~\cite[App.\ II]{gazi:cfht})
and because colimits of filtered diagrams are exact the sequence is exact.
Since tensor products commute with colimits, it is in fact pure exact and 
hence can be used to compute $\PExt^*(M,-)$.

As promised in the previous paragraph, we now show that every
module $M$ is a filtered colimit of finitely presented modules.
To avoid set theoretic problems, fix a set of finitely presented
modules containing a representative from each isomorphism class.
Let $\Lambda(M)$ be the
category whose objects are maps $P \ra M$ where $P$
is in our set of finitely presented modules.  
The morphisms are the obvious commutative triangles.
This category is filtered, and there is a natural functor 
$\Lambda(M) \ra \RMod$ sending $P \ra M$ to $P$.
The colimit of this diagram is $M$.
A smaller but less canonical filtered diagram of finitely presented
modules with colimit $M$ is described in~\cite[Exercise~I.2.10]{bo:ca}.

The exact sequence~\eqref{eq:ppr} leads to a spectral sequence
\begin{equation}\label{eq:ext-ss} 
E_2^{p,q} = \invlim^{\!p} \Ext^q(M_\alpha,N) \implies \Ext^{p+q}(M,N) 
\end{equation}
involving the derived functors of the inverse limit functor.
One way to construct this spectral sequence is as follows.
Break the exact sequence displayed above into short exact
sequences, defining modules $M_{i}$ in the process:
\[ 
  \xymatrixcolsep{0.6pc} 
  \xymatrixrowsep{1.05pc}
  \hspace*{-1em}
  \begin{array}{c}
  \xymatrix{ \relax 
   \makebox[2em][r]{$M = M_0$}
  && M_1 \ar@{ >->}[ld] && M_2 \ar@{ >->}[ld] && M_3 \ar@{ >->}[ld] \\
   & \Oplus_{\alpha} M_{\alpha} \ar@{{}->>}[lu] 
  && \Oplus_{\alpha \ra \beta} M_{\alpha} \ar@{{}->>}[lu] 
  && \Oplus_{\alpha \ra \beta \ra \gamma} M_{\alpha} \ar@{{}->>}[lu] 
  }
  \end{array}
  \cdots .
\]
Applying $\Ext^{*}(-,N)$ produces an unraveled exact couple
\[
  \hspace*{-1em}
  \begin{array}{c}
  \xymatrix@C-4.2pc@R-0.5pc{ \relax  
     \Ext^*(M,N) \ar[rd] && \Ext^*(M_1,N) \ar[ll] \circar[ll] \ar[rd] && \Ext^*(M_2,N) \ar[ll] \circar[ll] \ar[rd] && \Ext^*(M_3,N) \ar[ll] \circar[ll] \\
   & \Prod_{\alpha} \Ext^*(M_\alpha,N) \ar[ru]  
  && \Prod_{\alpha \ra \beta} \Ext^*(M_\alpha,N) \ar[ru]  
  && \Prod_{\alpha \ra \beta \ra \gamma} \Ext^*(M_\alpha,N) \ar[ru] 
  }
  \end{array}
  \!\!\cdots ,
\]
in which the horizontal maps are the connecting maps in the long
exact sequence of $\Ext$ groups.
This exact couple leads to a spectral sequence abutting to 
$\Ext^{*}(M,N)$, and
the $E_{2}$-term is the cohomology of the bottom row which,
by~\cite[App.\ II]{gazi:cfht}, is $\invlim^{\!p} \Ext^{q}(M_{\alpha}, N)$.
The same construction works if the sequence~\eqref{eq:ppr} is replaced by
any pure projective resolution of $M$.  The spectral sequences produced
in this way agree from the $E_{2}$-term onwards.
The $E_{2}$-term consists of the derived functors of $\Ext^{*}(-,N)$
with respect to the pure projective class.

The spectral sequence determines a decreasing filtration of $\Ext^{k}(M,N)$.
We write $\PlExt^{k}(M,N)$ for the $l$th stage, so
$\PzExt^{k}(M,N) = \Ext^{k}(M,N)$. 
The next stage, $\PoExt^{k}(M,N)$, consists of those extensions $\alpha$ 
in $\Ext^{k}(M,N)$ which can be factored into a product $\beta \gamma$ with
$\beta$ in $\Ext^{k-1}(K,N)$ and $\gamma$ in $\PExt^{1}(M,K)$ for some $K$.
Indeed, the map $\Ext^{k-1}(M_{1},N) \ra \Ext^{k}(M,N)$, whose image
is $\PoExt^{k}(M,N)$, is given by composition with the pure
extension $0 \ra M_{1} \ra \oplus_{\alpha} M_{\alpha } \ra M \ra 0$.
In general, for $0 \leq l \leq k$, $\PE{l}^{k}(M,N)$ consists of 
those extensions $\alpha$ in $\Ext^{k}(M,N)$
which can be factored into a product $\beta \gamma$ with
$\beta$ in $\Ext^{k-l}(K,N)$ and $\gamma$ in $\PExt^{l}(M,K)$ for some $K$.
Note that $\PkExt^{k}(M,N)$ is exactly the image of $\PExt^{k}(M,N)$
in $\Ext^{k}(M,N)$, and that 
$\PlExt^{k}(M,N)$ is zero for $l > k$.

\begin{question}
  From the exact couple it is clear that $\PoExt^{k}(M,N)$ can
  also be described as those extensions which pullback to zero
  under any map from a finitely presented module to $M$.
  Is it true that $\PlExt^{k}(M,N)$ consists of those
  extensions which pullback to zero under any map from a
  module of pure projective dimension at most $l-1$ to $M$?
\end{question}

\subsection{The relation between phantom maps and pure extensions}

Recall that $R$ is an associative ring and that $\D$ denotes the
derived category of (left) $R$-modules.
We begin by characterizing the finite objects in $\D$, \ie those
objects $X$ such that $\D(X,\oplus Y_{\alpha}) = \oplus \D(X,Y_{\alpha})$.
We first note that $R$, regarded as a complex concentrated in degree $0$,
is finite.
Indeed, $R$ is cofibrant, and so $\D(R,X) \iso \K(R,X)$.
It is easy to see that $\K(R,X) \iso H_{0}X$,
and since $H_{0}(\oplus Y_{\alpha }) = \oplus H_{0} Y_{\alpha }$, we see
that $R$ is finite, as claimed.
This also shows that $R$ is a weak graded generator for $\D$, \ie that
a complex $X$ is (isomorphic to) zero if and only if $\D(\Sigma^{n} R, X) = 0$ for all $n$.

The following concept will be of use to us.
An $R$-module $M$ is \dfn{FP} 
if there is a finite resolution of $M$ by finitely generated projectives.

We will also need the following terminology.
A full subcategory $\cT$ of a triangulated category $\cS$ 
is a \dfn{thick subcategory} if it is closed under cofibres,
retracts and desuspensions.
The thick subcategory \dfn{generated} by a collection $\cU$ of objects
is the full subcategory determined by all objects which can be
built from a finite number of objects of $\cU$ using cofibres,
retracts and desuspensions.
It is the smallest thick subcategory containing $\cU$.

\begin{prop}\label{pr:finites}
Let $X$ be an object of $\D$.  Then the following are equivalent:
\begin{roenumerate} %
\item $X$ is finite.
\item $X$ is in the thick subcategory generated by $R$.
\item $X$ is isomorphic to a bounded complex of finitely generated projectives.
\end{roenumerate}
Moreover, if $H_{n}X$ is FP for all $n$ and is zero for all but 
finitely many $n$, then $X$ is finite.
\end{prop}

\begin{note}
The converse of the last statement doesn't hold in general.
For example, if $R = \kay[x]/x^{2}$ for some field $\kay$,
then the chain complex 
\[
\cdots \lra 0 \lra \kay[x]/x^{2} \llra{x} \kay[x]/x^{2} \lra 0 \lra \cdots
\]
is finite, but has homology modules of infinite projective dimension.
\end{note}

\begin{proof}[Proof of Proposition~\ref{pr:finites}]
Since $R$ is a finite weak graded generator of $\D$,
it follows from~\cite[Corollary~2.3.12]{hopast:ash}
that (\emph{i}) and (\emph{ii}) are equivalent.
(That (\emph{ii}) implies (\emph{i}) is straightforward, but the other direction
is less so.)

We next prove that (\emph{iii}) implies (\emph{ii}): 
The thick subcategory generated by $R$ is closed
under finite coproducts and retracts, and so it contains all finitely
generated projective modules (considered as complexes concentrated
in one degree).
A bounded complex of finitely generated projectives can be built
from such complexes using a finite number of cofibres and thus
is also contained in the thick subcategory generated by $R$.

Now we prove that (\emph{ii}) implies (\emph{iii}):  
Let $\cT$ be the collection of complexes
isomorphic to a bounded complex of finitely generated projectives.
Since $R$ is in $\cT$, it suffices to prove that $\cT$ is a thick subcategory.
Given objects $X$ and $Y$ in $\cT$ and a map $f : X \ra Y$ we must show
that the cofibre of $f$ is in $\cT$.  We can assume without loss
of generality that $X$ and $Y$ are bounded complexes of finitely 
generated projectives.  A choice of cofibre has as its $n$th module
the direct sum $Y_{n} \oplus X_{n-1}$.  Thus the cofibre is again
in $\cT$.  The subcategory $\cT$ is clearly closed under desuspension, 
so it remains to show that $\cT$ is closed under retracts.
This is proved as Proposition~3.4 of~\cite{bone:hltc}.

Finally, we prove that if $X$ is a complex such that each
$H_{n}X$ is FP and only finitely many are non-zero, then
$X$ is finite.
If $X$ has no homology, then $X \iso 0$ in $\D$, so $X$ is finite.
If $X$ has homology concentrated in one degree, and this
module has a finite resolution by finitely generated projectives,
then $X$ is isomorphic to this finite resolution and is thus
finite (since (\emph{iii}) implies (\emph{i})).
Assume now that $X$ has non-zero homology only in a range
of $k$ degrees, with $k > 1$.
Without loss of generality, assume that $H_{n}X = 0$ for
\mbox{$n < 0$} and $n \geq k$.
Choose a finitely generated projective $P$ and a map $P \ra X$ inducing
an epimorphism $P \ra H_{0}X$ with FP kernel $K$.
Let $X'$ be the cofibre of the map $P \ra X$.
Then one finds that $H_{n}X'$ is zero for $n \leq 0$ and $n \geq k$.
Moreover, we have that $H_{n}X' = H_{n}X$ for $1 < n < k$ and
that there is a short exact sequence $0 \ra H_{1}X \ra H_{1}X' \ra K \ra 0$.
Since $H_{1}X$ and $K$ are FP, so is $H_{1}X'$,
and so by induction we can conclude that $X'$ is finite. 
The complex $P$ is certainly finite; therefore $X$ is finite as well
(since (\emph{iii}) implies (\emph{i})).
\end{proof}

\begin{note}\label{no:coherent}
If $R$ is a coherent ring over which every finitely presented
module has finite projective dimension, then the above result
simplifies.  First, over such a ring, a module is 
finitely presented if and only if it is FP.
Second, over a coherent ring, finitely presented modules 
form an abelian subcategory of the category of all modules,
and this subcategory is closed under retracts and extensions.
This allows one to show that a complex $X$ is finite if
and \emph{only if} each $H_{n}X$ is finitely presented and only
finitely many are non-zero.
In addition, in this situation,
the reliance on the result of B\"okstedt and Neeman
can be removed from the above proof.
(See~\cite[Exercise~I.2.11]{bo:ca} for a brief discussion
of coherence.  Note that Noetherian rings are coherent.)
\end{note}

Proposition~\ref{pr:finites} implies that there is a set of isomorphism
classes of finite objects and therefore that $\D$
has a phantom projective class (see Definition~\ref{de:gppc}
and the subsequent discussion).
The class $\Proj$ of projectives consists of all retracts of 
coproducts of finite objects,
and we write $\I$ for the ideal of phantom maps.
Since $R$ is finite, the phantom projective class generates.

We recall a standard fact which is easily proved.

\begin{prop}\label{pr:ext}
  Let $M$ and $N$ be $R$-modules.  Then
  \[ \D(M,\Sigma^k N) \iso \Ext_R^k(M,N) . \]
\\*[-.52in]
\qed
\end{prop}

We can now prove one of the main results of this section.

\begin{thm}\label{th:agree}
Let $M$ be a filtered colimit of FP modules and let $N$ be any $R$-module.
Then the phantom spectral sequence abutting to $\D(M,N)$ 
is the same as the spectral sequence~\eqref{eq:ext-ss} described in
the previous part of this section.
In particular, the filtrations agree:
\[ \I^l(M,\Sigma^k N) \iso \PlExt^k(M,N) . \]
\end{thm}

When we say that the spectral sequences are the same, we mean
that they are naturally isomorphic from the $E_{2}$-term onwards.

\begin{proof}
Let $\{M_{\alpha}\}$ be a filtered diagram of FP modules
with a colimiting cone to $M$.
Then, regarding these modules as complexes concentrated in degree zero,
the cone is a minimal cone.
Indeed, the cone from the complexes $\{M_{\alpha}\}$ to the complex $M$ 
becomes a colimiting cone under $\D(\Sigma^{n} R,-) = H_{n}(-)$
for each $n$.
And since filtered colimits are exact, one can use the five-lemma
to show that it becomes a colimiting cone under $\D(W,-)$ for
each finite $W$.
This is what it means for the cone to be a minimal cone.

By Proposition~\ref{pr:finites}, the complexes $M_{\alpha}$ are finite.

We saw in Theorem~\ref{th:limi} that from a minimal cone on a 
filtered diagram of finite objects one can construct a phantom
resolution.  In fact, the construction corresponds exactly to
the construction of a pure resolution of $M$ in the previous
part of this section.
Moreover, to get the spectral sequence~\eqref{eq:ext-ss} we apply
the functor $\Ext^{*}(-,N)$.
To get the phantom spectral sequence we apply the functor
$\D(-,N)_{*}$.  By Proposition~\ref{pr:ext}, these agree.
Thus the spectral sequences agree.  
\end{proof}

The following lemma will allow us to find modules of pure projective
dimension greater than one.

\begin{lemma}\label{le:flat}
Let $M$ be a flat $R$-module.  
Then any projective resolution of $M$ is a pure projective resolution.
Moreover, the projective dimension of $M$ equals the pure projective dimension
of $M$ and for any $N$ the natural map $\PExt^*(M,N) \ra \Ext^*(M,N)$ is
an isomorphism.
\end{lemma}

\begin{proof}
A projective resolution of $M$ is an exact sequence
\[
0 \lla M \lla P_{0} \lla P_{1} \lla \cdots
\]
with each $P_{i}$ projective (and hence pure projective).
We must show that this sequence is pure exact.
Write $M_{i}$ for the image of the map $P_{i-1} \la P_{i}$.
For each right module $E$,
the sequence $0 \la M \la P_{0} \la M_{1} \la 0$ is exact
under $E \tensor -$ since $M$ is flat.
And since $M$ and $P_{0}$ are flat, so is $M_{1}$.
Thus we can inductively conclude that each sequence
$M_{i} \la P_{i} \la M_{i+1}$ is pure exact.
Therefore, the resolution is pure exact and we have proved
the first part of the proposition.

From what we have proved, it follows that a projective resolution
of $M$ can be used to calculate both $\PExt^*(M,N)$ and $\Ext^*(M,N)$.
Thus the natural map $\PExt^*(M,N) \ra \Ext^*(M,N)$ is an isomorphism,
and the projective and pure projective dimensions of $M$ agree.
\end{proof}

For the following theorem, we need to define the notion of a 
``pure complex''.

\begin{defn}
A complex $X$ is \dfn{pure} if for each $n$ the short exact
sequences
\[ 0 \lra B_{n}X \lra Z_{n}X \lra H_{n}X \lra 0 \]
and
\[ 0 \lra Z_{n}X \lra X_{n} \lra B_{n-1}X \lra 0 \]
are pure exact.
\end{defn}

\begin{thm}\label{th:compare}
Let $R$ be a coherent ring over which finitely presented modules 
have finite projective dimension.
\begin{roenumerate}
\item If every $R$-module has pure projective dimension less than $n$,
then every pure complex has length at most $n$.
\item If every pure complex has length at most $n$, then 
every \emph{flat} $R$-module has pure projective dimension less than $n$.
\item If every pure complex has phantom projective dimension at most one,
then every $R$-module has pure projective dimension at most one.
\end{roenumerate}
\end{thm}

\begin{proof}
The proof of (\emph{i}) is similar to the proof of Theorem~\ref{th:algh}.
Note, however, that we do not claim that each pure complex has
phantom projective dimension less than $n$.
The correct statement along these lines must be made in the
``pure'' derived category, and will be described in a subsequent paper.

We proceed to prove (\emph{ii}).
Assume that there exists a flat module $M$
of projective dimension at least $n$.
(Recall that for flat modules, the pure projective and projective
dimensions agree.)
Then $\Ext^{n}(M,N)$ is non-zero for some $N$.
But $\I^{n}(M,\Sigma^{n}N) = \PE{n}^{n}(M,N) = \Ext^{n}(M,N)$.
Indeed, the first equality is Theorem~\ref{th:agree} and
the second is Lemma~\ref{le:flat}.
(Note that the hypothesis of Theorem~\ref{th:agree} is satisfied
because we have assumed that every finitely presented $R$-module
is FP.)
Thus there are $n$ phantoms in $\D$ which compose non-trivially,
and the source is $M$, a pure complex.

Finally, we prove (\emph{iii}).
Let $M$ be an $R$-module.
By assumption there is a phantom exact sequence
\[ 0 \lra P \lra Q \lra M \lra 0 \]
with $P$ and $Q$ retracts of sums of finite complexes.
By Note~\ref{no:coherent}, $H_{0}P$ and $H_{0}Q$ are finitely
presented.  And since a finitely presented module is finite
when regarded as a complex concentrated in degree zero, it
is easy to see that 
\[ 0 \lra H_{0}P \lra H_{0}Q \lra M \lra 0 \]
is a pure projective resolution of $M$.
Thus, in this setting, every $R$-module has pure projective dimension
at most one.
\end{proof}

\begin{cor}
If $R$ is a coherent ring over which finitely presented modules are FP
and whose derived category $\D$ is a Brown category
(see Definition~\ref{de:br}), 
then every $R$-module has pure projective dimension at most one.
\end{cor}

\begin{proof}
Suppose that $\D$ is a Brown category.
Then, by Theorem~\ref{th:pdim1}, it follows that every
complex has phantom projective dimension at most one.
Thus, by Theorem~\ref{th:compare}~(\emph{iii}), every
$R$-module has pure projective dimension at most one.
\end{proof}

It is proved in~\cite[Theorem~4.1.5]{hopast:ash}
and in~\cite[Section~5]{ne:tba} that 
when $R$ is a countable ring, $\D$ is a Brown category.
So we have shown in particular that a countable coherent ring over
which finitely presented modules are FP has pure global dimension
at most one.  In fact, \emph{any} countable ring has pure global
dimension at most one~\cite[Proposition~10.5]{grje:dcrf}.

\begin{example}
Let $\kay$ be an uncountable field.
Write $\kay[x,y]$ for the polynomial ring and $\kay(x,y)$ for its
field of fractions.
A theorem of Kaplansky~\cite{ka:hdqf} states that the projective
dimension of $\kay(x,y)$ as a $\kay[x,y]$-module is at least two.
Moreover, $\kay(x,y)$ is a flat $\kay[x,y]$-module and $\kay[x,y]$
is a coherent ring of global dimension 2.
Thus, by Theorem~\ref{th:compare}~(\emph{ii}), 
in the derived category of $\kay[x,y]$-modules, $\kay(x,y)$
is the source of a non-zero composite of two phantoms.
In particular, the derived category of $\kay[x,y]$-modules is
not a Brown category.
I learned this example from Neeman~\cite{ne:tba} and Neeman credits it to
Bernhard Keller.

The phantom maps which compose non-trivially when there is a flat
module of projective dimension greater than one can be made more explicit.
Let $M$ be a flat $R$-module of projective dimension at least $n$
and let
\[
  0 \lra N \lra A_1 \lra \cdots \lra A_n \lra M \lra 0 
\]
be a pure exact sequence which is non-zero as an element of $\Ext^{n}(M,N)$.
(This is possible by Lemma~\ref{le:flat}.)
This sequence factors into $n$ short exact sequences which are also
pure exact.  By Theorem~\ref{th:agree} each of these represents a phantom map;
their composite is the given non-zero element of 
$\Ext^{n}(M,N) = \D(M,\Sigma^{n}N)$.
\end{example}

\renewcommand{\baselinestretch}{1}\normalsize

\newpage

\noindent
\textbf{\large This version has been updated compared to the published version:}

\medskip

\noindent
The published version appeared in {\em Advances in Mathematics} 136 (1998), 284-339.
Since then, the following changes have been made:

\begin{itemize}
\item In Proposition~\ref{pr:cones}~(ii), the assumption that the
      objects in the diagram are projective was added, and two 
      sentences were added at the end of the proof using this assumption.
      The discussion before Theorem~\ref{th:milnor} is adjusted accordingly.
\item In Section~\ref{se:ass}, the assumption that the ideal $\I$ is closed
      under countable coproducts was added to Propositions~\ref{pr:condconv}
      and~\ref{pr:collapse}.  And in Section~\ref{sse:ucss}, the assumption
      that the forgetful functor $U$ commutes with coproducts was added.
      The introduction was modified accordingly.
      Thanks to Ciprian Modoi and Ralf Meyer for pointing out this problem.
\item My address and e-mail address have been updated.
\item Reference \cite{bo:ccss} has been updated.
\item Reference \cite{br:rhaass}, which I was unaware of and should have
      cited, has been added.
\item Various minor typos fixed.
\end{itemize}

\end{document}